\documentclass{article}
\usepackage{amssymb}
\usepackage{mathrsfs}
\usepackage{cite}

\usepackage{stmaryrd}
\usepackage{amsfonts}
\usepackage{amsmath}
\usepackage{bm}
\usepackage{indentfirst}
\usepackage{array}
\usepackage{arydshln}
\usepackage{graphicx}
\usepackage{float}
\usepackage{subfigure}
\usepackage{listings}
\usepackage{epstopdf}
\numberwithin{equation}{section}
\usepackage{geometry}
\usepackage{amsthm}
\usepackage{mathrsfs}
\usepackage[verbose]{hyperref}
\newtheorem{thm}{Theorem}

\newtheorem{lem}{Lemma}

\newtheorem{coro}{Corollary}
\makeatletter
\newcommand{\rmnum}[1]{\romannumeral #1}
\newcommand{\Rmnum}[1]{\expandafter\@slowromancap\romannumeral #1@}
\makeatother
\title{\Large \bf{Hopf and double Hopf bifurcations in a delayed lateral vibration model of footbridges induced by pedestrians
\footnote{This work is supported by the NNSF (12471158) of China (Xuemei Li, Email:lixuemei\underline{~}1@sina.com) and the NSF (2023JJ40659) of Hunan Province, China (Yechi Liu, Email:lyc9009@sina.cn).}}
\setcounter{footnote}{-1}
\author{Xuemei Li,\qquad Yechi Liu\\
{\small Key Laboratory of High Performance Computing and Stochastic Information Processing}\\
{\small Department of Mathematics, Hunan Normal University, Changsha, Hunan 410081, China}\\
}}
\date{}
\begin{document}\large
\maketitle

\textbf{Abstract}. In this paper, we investigate the dynamical behaviors of a delayed lateral vibration model of footbridges proposed based on the facts that pedestrians will reduce their walking speed or stop walking when the response of the footbridge becomes sufficiently large, and that the bridge velocity can not be changed at once when the pedestrians begin to walk on the bridge. By analyzing the distribution of roots of the associated characteristic equation, we find that there are only two types of bifurcations in this model: Hopf bifurcation and double Hopf bifurcation, and give the condition on the stability of the trivial solution. By using the center manifold theorem and bifurcation theory of delayed differential equations, we obtain the dynamical behavior in these bifurcations, specially including the stability of periodic solutions and invariant tori bifurcating from the trivial solution in these bifurcations. Finally, we prove that this model exhibits quasi-periodic vibrations by KAM theorems, besides periodic vibrations.\\

\textbf{Keywords}. Lateral vibration model, Hopf bifurcation, Double Hopf bifurcation, Quasi-periodic vibration.

\section{Introduction}
Walking human-induced vibration on floors, such as footbridges, stadia, long-span floors, etc., has been studied since Harper \cite{Har62} carried out a walking test in 1962 and proposed a M-shape force model. Andriacchi et al. \cite{AOG77} in 1977 proposed a three-dimensional walking model based on measuring the footfall forces in the three orthogonal directions of vertical, longitudinal and lateral components. Because of the use of new materials and constructure technology, some footbridges are becoming lighter, longer and more slender, and consequently more susceptible to vibration. Walking human-induced vibration on footbridges has drawn considerable attention. By using different measuring devices, the walking models on footbridges were studied in \cite{Wh82,RPA88,KB01,EHSP96,SPR16,CLCH23}. In \cite{MP21,BBR16,SP18,KKRE18,RB11}, some methods of determining the vertical component force of walking and mathematical models are proposed. For the longitudinal walking force, Kumar et al. \cite{KKR18} gave a model based on self-sustained oscillator. For more studies on longitudinal component, see literature reviews \cite{ZPR05}.

These researches on vibrations induced by pedestrians before 2000 had mainly focused on vertical vibrations of footbridges which gave some suggestions in structural design and calculation of bridges. The problem of the excessive lateral vibration of bridges induced by crowds of pedestrians has caught many researchers' attention and triggered a wave of study on the human-structure interaction in the lateral direction since the Paris Softerino Bridge and the London Millennium Bridge experienced large amplitude lateral vibrations on their opening days \cite{RD16,IGJ12,ZPR05,RPB09,FD16}. On one hand, there have been many efforts to measure forces due to walking experiments in order to understand the footbridge-pedestrian foot interaction \cite{RPB09,ZPR05,IGRJ11,RP07,BCB18}. On the other hand, according to the experimental studies for human-induced forces, many lateral walking force models in biomechanics have been proposed to explore the mechanics of human balance during walking and energy exchange between a crowd and a floor, such as Mcdonald's et al. proposed the inverted pendulum model \cite{Ma08,BMB13,BMB12,Gol15} based on mathematical pendulum equation, Belykh et al. \cite{BJB17,BDB22} investigated pedestrian-induced bridge instability (i.e. phase locking and bridge wobbling appear) based on the inverted pendulum model and a van der Pol-type oscillator model; Kumar et al. \cite{KKE16} proposed the modified hybrid van der Pol-Duffing-Rayleigh oscillator model based on classifications in \cite{ETA13}, where the footbridge-Pedestrian interaction is modelled as uncoupled problem, and proved the existence of a stable limit cycle of their model by using the energy balance method and Lindstedt-Poincare perturbation technique. In \cite{HZJZ21}, the lateral vibration of a bridge subjected to walking pedestrians is modelled by a partial differential equation based on considering the bridge to be a slender beam and the theory of coupled oscillators.

Based on the analysis in \cite{DFF01} and the fact that pedestrians will reduce their walking speed or stop walking when the response of the footbridge becomes sufficiently large, Nakamura \cite{Nak04,NK09} proposed the following model
\begin{equation}\label{11}
M_f\ddot x+C_f\dot x+K_fx=F_p(t),
\end{equation}
where $M_f,C_f$ and $K_f$ are modal mass, damping and stiffness of a footbridge in the lateral direction, respectively; $x,\dot x$ and $\ddot x$ are modal lateral displacement, velocity and acceleration of the footbridge, respectively; $F_p(t)$ is the modal lateral force induced by pedestrians on the footbridge given by
\begin{equation}\label{12}
F_p(t)=K_1K_2M_pgG(f)H_p(\dot x),
\end{equation}
where $K_1$ is the ratio of the lateral force to the pedestrians' weight, $K_2$ is the percentage of pedestrians who synchronize to the lateral bridge vibration, $M_pg$ is the modal self-weight of pedestrians, $G(f)$ is the function to describe how pedestrians synchronize with the bridge natural frequency and is assumed $G(f)=1$, $H_p(\dot x)=\dot x/(K_3+|\dot x|)$ is a function to describe the pedestrians' synchronization nature, which is assumed to be proportional to the footbridge velocity $\dot x$ at low velocities, $K_3$ is used to describe the saturation rate. In \cite{ZXX13}, $H(\dot x)$ is taken the hyperbolic tangent function ${\rm tanh}(K_3\dot x)$, which satisfies all the assumptions about $H_p(\dot x)$ defined in \eqref{12}, that is, Zhen et al. \cite{ZXX13} considered the lateral vibration model (\ref{11}-\ref{12}) of footbridges induced by pedestrians given by
\begin{equation}\label{13}
M_f\ddot x+C_f\dot x+K_fx=K_1K_2M_pgG(f){\rm tanh}(K_3\dot x).
\end{equation}
The existence of periodic solutions to \eqref{13} and the analytical approximate expression are investigated in \cite{ZXX13,ZXS16} by using the Hopf bifurcation theory and the energy method.

Considering the facts that the bridge velocity can not be changed at once when the pedestrians begin to walk on the bridge, and that the prediction results Nakamura's model were significantly affected by the time delay, Zhen et al. \cite{ZXX13} introduced the time delay in Nakamura's model, which leads (\ref{11}-\ref{12}) becoming
\begin{equation}\label{14}
\ddot x(t)+\alpha_2\dot x(t)+\alpha_1x(t)=\kappa{\rm tanh}\big(K_3\dot x(t-\tau)\big),
\end{equation}
where
\begin{equation*}
\alpha_1=\frac{K_f}{M_f}>0,\qquad \alpha_2=\frac{C_f}{M_f}>0,\qquad \Psi=\frac{K_1K_2M_pg}{M_f}>0,
\end{equation*}
$\tau$ is the time delay of the interaction between pedestrians and the footbridge. They performed the qualitative analysis on the delayed model \eqref{14}. Subsequently, Zhen, et al. \cite{ZXS16} performed the quantitative analysis on the existence of periodic solutions of \eqref{14} by the energy method, as well as the influence of the time delay on lateral vibration amplitude of the footbridge (i.e. periodic solutions). By regarding the time delay $\tau$ as the bifurcation parameter, and obtaining the critical value $\tau_0$ and the transversality condition \big(i.e. $\left.\frac{{\rm d}}{{\rm d}\tau}({\rm Re}\lambda)\right|_{\tau=\tau_0}\neq0$, $\lambda$ is the eigenvalue of \eqref{14}\big) corresponding to the Hopf bifurcation, they gave the existence of the bifurcating limit cycle (i.e. the periodic solution) and under the assumption that the periodic solution of \eqref{14} has the form $x(t)=r\cos(\omega t+\phi)$, calculated a first-order approximation of the periodic solution by using the energy method. In this paper, we will rigorously prove the existence and stability of limit cycles bifurcating from Hopf bifurcation in \eqref{14}, and obtain approximate expressions of the corresponding periodic solutions.

The aim of this paper is to investigate bifurcation phenomena in \eqref{14}. By analyzing the distribution of eigenvalues of \eqref{14}, we find that there are only two types of bifurcations in \eqref{14}: Hopf bifurcation and double Hopf bifurcation, and there exist many parameter bifurcation values (i.e. critical values) at which one of these two bifurcations occurs. By using the center manifold theorem and bifurcation theory of delayed differential equations (DDEs), we will analyze the dynamical behavior in these bifurcations. In particular, we will prove that \eqref{14} exhibits abundant periodic and quasi-periodic vibrations. But \eqref{13}, without the time delay, only has periodic vibration. Our results show that the time delay can cause quasi-periodic vibration.

Obviously, \eqref{14} is a class of delayed oscillator models. There are some relevant works on the bifurcation analysis for oscillator models with the time delay. For the van der Pol's oscillator with a nonlinear delay feedback
\begin{equation}\label{15}
\ddot x(t)+c\big(x^2(t)-1\big)\dot x(t)+x(t)=c f\big(x(t-\tau)\big),
\end{equation}
where $f\in C^3,f(0)=f^{\prime\prime}(0)=0,f^\prime(0)=a\neq0$, Jiang and Wei \cite{JW08,WJ05} discussed the fold and Hopf bifurcations. Subsequently, more bifurcation phenomena in \eqref{15} are obtained, including Bogdanov-Takens bifurcation \cite{JY07}, zero-Hopf bifurcation \cite{WJ10,WW11}, triple-zero bifurcation \cite{HLS12}, non-semisimple 1:1 resonance Hopf bifurcation \cite{ZG13}. Recently, Li and Yu \cite{LY24} analyzed the double Hopf bifurcation, and obtained the quasi-periodic oscillation phenomenon in \eqref{15} by using KAM theory. Brambuger et al. \cite{BDL14} studied the zero-Hopf bifurcation in a more general model with a delay feedback forcing
\begin{equation}\label{16}
\ddot x(t)+c\big(x^2(t)-1\big)\dot x(t)+x(t)=g\big(\dot x(t-\tau),x(t-\tau)\big),
\end{equation}
where $g\in C^3$,
\begin{equation}\label{g}
g(0,0)=0,\quad g_{\dot x}(0,0)=a\neq0,\quad g_x(0,0)=b\neq0
\end{equation}
(with the notion $g_y=\frac{\partial g}{\partial y}$). For hybrid models of van der Pol and Duffing oscillators with the time delay, there are also similar results on bifurcation phenomena, see \cite{DJY13,JZS15,MLF08}, for example.

The rest of this paper is organized as follows. In Section 2, we analyze the characteristic equation of \eqref{14} and obtain the distribution of the eigenvalues, which implies that there are only two types of eigenvalues with zero real parts: a pair of simple purely imaginary eigenvalues and two pairs of simple purely imaginary eigenvalues, which maybe induce Hopf bifurcation and double Hopf bifurcation in \eqref{14}, respectively. We also investigate the stability of the zero solution of \eqref{14}. In the second half of this section, we give a decomposition of delayed differential equations with parameters which is used to analyze Hopf and double Hopf bifurcations in \eqref{14}. In Section 3, we discuss the direction and stability of the bifurcating periodic solutions from Hopf bifurcation regarding the time delay $\tau$ and coefficients of the linear part in \eqref{14} as the bifurcation parameters, respectively. In Section 4, we study double Hopf bifurcation in \eqref{14} and obtain the bifurcation diagram of the truncated normal form corresponding to \eqref{14}, and specially get a parameter region where the truncated system has a quasi-periodic solution. Then, by using KAM theorems, we prove that \eqref{14} also has a quasi-periodic solution for most parameter values of this region. Appendix collects some expressions coming from calculating normal forms.

\section{Stability and bifurcation values of \eqref{14}}
In this section, we obtain possible bifurcation values and the stability of the zero solution to \eqref{14} by analyzing the associated characteristic equation. We also give a decomposition of general delayed differential equations with parameters near a bifurcation value on account of the center subspace.

\subsection{The distribution of eigenvalues and stability}
Without loss of generality, we take $\alpha_1=1$ in \eqref{14}. Otherwise, we make a transformation
\begin{equation*}
s=\sqrt{\alpha_1}t,\qquad y(s)=x(\frac{s}{\sqrt{\alpha_1}})
\end{equation*}
and rewrite $y,\frac{\alpha_2}{\sqrt{\alpha_1}},\frac{\kappa}{\alpha_1},\sqrt{\alpha_1}k_3$ and $\sqrt{\alpha_1}\tau$ as $x,\alpha_2,\kappa,k_3$ and $\tau$, respectively, and \eqref{14} can be transformed into
\begin{equation}\label{21}
\ddot x(t)+\alpha_2\dot x(t)+x(t)=\kappa{\rm tanh}\big(K_3\dot x(t-\tau)\big).
\end{equation}
Letting $x(t)=x_1(t),\dot x(t)=x_2(t)$. \eqref{21} can be written as
\begin{equation}\label{22}
\left\{\begin{aligned}
&\dot{x_1}(t)=x_2(t),\\
&\dot{x_2}(t)=-x_1(t)-\alpha_2x_2(t)+\kappa{\rm tanh}\big(K_3x_2(t-\tau)\big).
\end{aligned}\right.
\end{equation}
Noting
\begin{equation*}
\kappa{\rm tanh}(K_3y)=\alpha_3y+\alpha_4y^3+\alpha_5y^5+O(y^7)
\end{equation*}
with
\begin{equation*}
\alpha_3=\kappa K_3>0,\alpha_4=-\frac{1}{3}\kappa K_3^3<0,\alpha_5=\frac{2}{15}\kappa K_3^5,
\end{equation*}
and the origin $(0,0)$ is an equilibrium, the associated characteristic equation is
\begin{equation}\label{23}
h(\lambda,\tau):=\lambda^2+\alpha_2\lambda+1-\alpha_3\lambda e^{-\lambda\tau}=0
\end{equation}
with $\alpha_2>0,\alpha_3>0$ and $\tau\geqslant0$. Obviously, $\lambda=0$ is not a root of \eqref{23}.

\vskip 0.1in

\noindent (a). For the case of $\tau=0$, the root $\lambda$ of \eqref{23} satisfies ${\rm Re}\lambda>0$ for $\alpha_3>\alpha_2$; ${\rm Re}\lambda=0$ for $\alpha_3=\alpha_2$; ${\rm Re}\lambda<0$ for $\alpha_3<\alpha_2$.

\vskip 0.1in

When $\tau=0$, letting $\lambda={\rm i}\omega$ with $\omega>0$ be a root of \eqref{23} with a zero real part, substituting it into \eqref{23} and separating the real and imaginary parts imply
\begin{equation}\label{24}
\left\{\begin{aligned}
&1-\omega^2=\alpha_3\omega\sin(\omega\tau),\\
&\alpha_2=\alpha_3\cos(\omega\tau).
\end{aligned}\right.
\end{equation}
It follows from \eqref{24} that $\omega$ satisfies
\begin{equation}\label{25}
\omega^4+(\alpha_2^2-\alpha_3^2-2)\omega^2+1=0.
\end{equation}
It is easy to see that

\vskip 0.1in

\noindent (b). If $\alpha_3>\alpha_2$, \eqref{25} has two positive roots
\begin{equation}\label{26}
\omega_\pm=\frac{1}{2}(\sqrt{4+\alpha_3^2-\alpha_2^2}\pm\sqrt{\alpha_3^2-\alpha_2^2});
\end{equation}

\vskip 0.1in

\noindent (c). If $\alpha_3=\alpha_2$, \eqref{25} only has a positive root 1;

\vskip 0.1in

\noindent (d). If $\alpha_3<\alpha_2$, \eqref{25} does not have a positive root. Hence, for any $\tau\geqslant0$, \eqref{23} does not have a purely imaginary root if $\alpha_3<\alpha_2$.

\vskip 0.2in

As for $\alpha_3>\alpha_2$,
\begin{equation*}
1-\omega_\pm^2=\frac{1}{2}\left(\alpha_2^2-\alpha_3^2\mp\sqrt{(\alpha_3^2-\alpha_2^2)(4+\alpha_3^2-\alpha_2^2)}\right),\qquad 1-\omega_+^2<0,\quad 1-\omega_-^2>0,
\end{equation*}
\eqref{24} implies
\begin{equation*}
\omega_-\tau\in(2\pi j,2\pi j+\frac{\pi}{2}),\quad \omega_+\tau\in(2\pi (j+1)-\frac{\pi}{2},2\pi (j+1)),\qquad j=0,1,2,\cdots.
\end{equation*}
Thus, we obtain the following lemma for the bifurcation values of \eqref{22}.
\begin{lem}\label{lem1}
{\rm (\rmnum{1})}. If $\alpha_3>\alpha_2$, then there exists a sequence $\{\tau_j^\pm\}$ with $\tau_{j+1}^\pm>\tau_j^\pm$ such that \eqref{23} has a pair of simple purely imaginary roots $\pm{\rm i}\omega_\pm$ when $\tau=\tau_j^\pm$, respectively, where $\omega_\pm$ are defined by \eqref{26},
\begin{equation}\label{26+}
\tau_j^+=\frac{1}{\omega_+}(2\pi-\arccos\frac{\alpha_2}{\alpha_3}+2\pi j),\quad \tau_j^-=\frac{1}{\omega_-}(\arccos\frac{\alpha_2}{\alpha_3}+2\pi j),\qquad j=0,1,2,\cdots.
\end{equation}
{\rm (\rmnum{2})}. If $C_1:\,\alpha_3=\alpha_2$, then there exists a sequence $\{\tau_j\}$ with $\tau_{j+1}>\tau_j$ such that \eqref{23} has a pair of simple purely imaginary roots $\pm{\rm i}$ when $\tau=\tau_j$,  where $\tau_j=2\pi j,j=0,1,2,\cdots$.
\end{lem}
Let $\lambda(\tau)=\sigma(\tau)+{\rm i}\omega(\tau)$ be a root of \eqref{23}. Substituting $\lambda(\tau)$ into \eqref{23} and taking the deribative with respect to $\tau$, we have the following lemma.
\begin{lem}\label{lem2}
Let $\lambda(\tau)=\sigma(\tau)+{\rm i}\omega(\tau)$ be a root of \eqref{23} satisfying $\sigma(\tau_0)=0$ and $\omega(\tau_0)=\omega_+,\omega_-,1$ for $\tau_0=\tau_j^+,\tau_j^-,\tau_j(j=0,1,2,\cdots)$, respectively. Then
\begin{align*}
&\sigma^\prime(\tau_j^+)=\frac{(\omega_+^2-1)(\omega_+^2+1)}{4+\alpha_3^2-\alpha_2^2+2\tau_j^+\alpha_2+\tau_j^+(\tau_j^+\alpha_3^2+2\alpha_2)\omega_+^2}>0,\\
&\sigma^\prime(\tau_j^-)=\frac{(\omega_-^2-1)(\omega_-^2+1)}{4+\alpha_3^2-\alpha_2^2+2\tau_j^-\alpha_2+\tau_j^-(\tau_j^-\alpha_3^2+2\alpha_2)\omega_-^2}<0,\\
&\sigma^\prime(\tau_j)=0,\qquad \sigma^{\prime\prime}(\tau_j)=-\frac{4\alpha_2}{(1+\alpha_2\tau_j)(2+\alpha_2\tau_j)^2}<0,\\
&\omega^\prime(\tau_j^\pm)=-\frac{\omega_\pm(\tau_j^\pm\alpha_3^2\omega_\pm^2+\alpha_2\omega_\pm^2+\alpha_2)}{4+\alpha_3^2-\alpha_2^2+2\tau_j^\pm\alpha_2+\tau_j^\pm(\tau_j^\pm\alpha_3^2+2\alpha_2)\omega_\pm^2}<0,\\
&\omega^\prime(\tau_j)=-\frac{\alpha_2}{2+\tau_j\alpha_2}<0
\end{align*}
for $j=0,1,2,\cdots$.
\end{lem}

\vskip 0.2in

Next, we give the arrangement of these two sequences $\{\tau_j^+\}$ and $\{\tau_j^-\}$. Noting $\alpha_3>\alpha_2$,
\begin{equation*}
\tau_{j+1}^+-\tau_j^+=\frac{2\pi}{\omega_+}<\frac{2\pi}{\omega_-}=\tau_{j+1}^--\tau_j^-,\qquad j=0,1,2,\cdots
\end{equation*}
and
\begin{equation*}
\tau_0^+-\tau_0^-=\sqrt{a}\left(\pi\Big(1-\sqrt{\frac{b}{a}}\Big)-\arccos\frac{\alpha_2}{\alpha_3}\right),
\end{equation*}
where $a=4+\alpha_3^2-\alpha_2^2$ and $b=\alpha_3^2-\alpha_2^2$, we obtain the following lemma.
\begin{lem}\label{lem3}
Suppose that $\alpha_3>\alpha_2>0$ is satisfied.\\
{\rm (\rmnum{1})}. If
\begin{equation}\label{27}
\frac{\alpha_2}{\alpha_3}>\cos\left(\pi\Big(1-\sqrt{\frac{\alpha_3^2-\alpha_2^2}{4+\alpha_3^2-\alpha_2^2}}\Big)\right),
\end{equation}
then $\tau_0^+>\tau_0^-$ and there exists an integer $m\geqslant0$ such that
\begin{equation*}
0<\tau_0^-<\tau_0^+<\tau_1^-<\tau_1^+<\cdots<\tau_m^-<\tau_m^+<\tau_{m+1}^+\leqslant\tau_{m+1}^-
\end{equation*}
and $\tau_j^+<\tau_j^-$ for $j\geqslant m+2$.\\
{\rm (\rmnum{2})}. If
\begin{equation}\label{28}
C_2:\quad \frac{\alpha_2}{\alpha_3}=\cos\left(\pi\Big(1-\sqrt{\frac{\alpha_3^2-\alpha_2^2}{4+\alpha_3^2-\alpha_2^2}}\Big)\right),
\end{equation}
then $0<\tau_0^+=\tau_0^-$, and $\tau_j^+<\tau_j^-$ for $j\geqslant1$.\\
{\rm (\rmnum{3})}. If
\begin{equation}\label{29}
\frac{\alpha_2}{\alpha_3}<\cos\left(\pi\Big(1-\sqrt{\frac{\alpha_3^2-\alpha_2^2}{4+\alpha_3^2-\alpha_2^2}}\Big)\right),
\end{equation}
then $0<\tau_j^+<\tau_j^-$ for $j\geqslant0$.
\end{lem}
By conclusions (a) and (d), Lemmas \ref{lem1}-\ref{lem3} and Rouch\'{e}'s theorem, we obtain the following theorem on the stability of the zero solution of \eqref{22}.
\begin{thm}\label{thm1}
For the system \eqref{22},\\
{\rm (\rmnum{1})}. If $\alpha_3<\alpha_2$, then all roots of \eqref{23} have negative real parts and the zero solution of \eqref{22} is asymptotically stable for all $\tau\geqslant0$;\\
{\rm (\rmnum{2})}. If $\alpha_3=\alpha_2$, then all roots of \eqref{23} have negative real parts and the zero solution of \eqref{22} is asymptotically stable for $\tau\neq\tau_j,j=0,1,2,\cdots$;\\
{\rm (\rmnum{3})}. If $\alpha_3>\alpha_2$ and \eqref{27} holds true, then there is an integer $m\geqslant0$ such that when $\tau\in(\tau_j^-,\tau_j^+)$ with $j=0,1,\cdots,m$, all roots of \eqref{23} have negative real parts and the zero solution of \eqref{22} is asymptotically stable; when $\tau>\tau_m^\pm$ or $\tau\in(\tau_{j-1}^+,\tau_j^-)$ with $j=0,1,\cdots,m$ and $\tau_{-1}^+=0$, \eqref{23} has at least a pair of roots with positive real parts and the zero solution of \eqref{22} is unstable;\\
{\rm (\rmnum{4})}. If $\alpha_3>\alpha_2$ and \eqref{28} or \eqref{29} holds true, then \eqref{23} has at least a pair of roots with positive real parts and the zero solution of \eqref{22} is unstable for all $\tau\geqslant0$.
\end{thm}

\begin{figure}[H]
  \centering
  \includegraphics[width=2in]{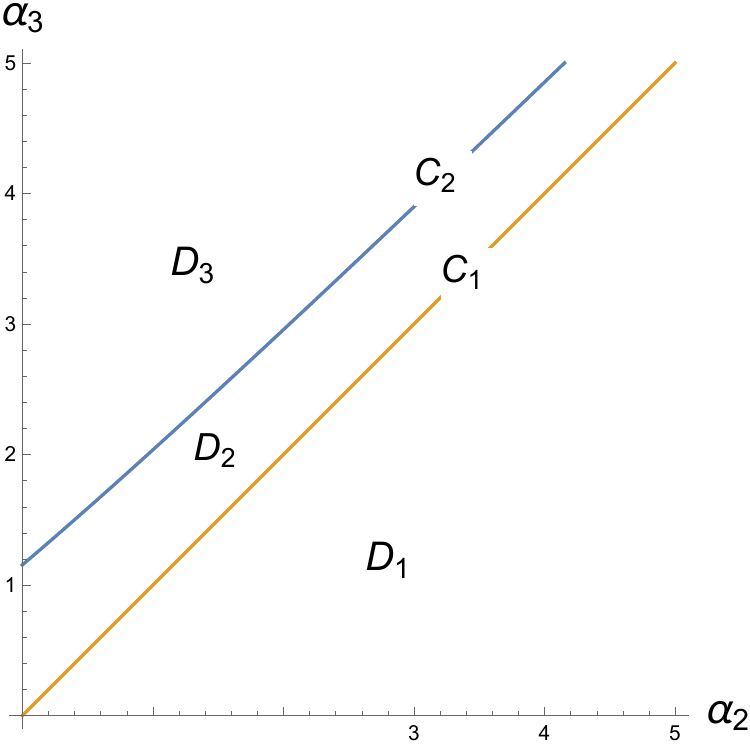}
  \caption{The decomposition of the $(\alpha_2,\alpha_3)$-parameter plane on the stability of the zero solution of \eqref{22}. The curves $C_1:\alpha_3=\alpha_2$ and $C_2:\frac{\alpha_2}{\alpha_3}=\cos\left(\pi(1-\sqrt{\frac{\alpha_3^2-\alpha_2^2}{4+\alpha_3^2-\alpha_2^2}})\right)$ divide the first quadrant of the $(\alpha_2,\alpha_3)$-parameter plane into three regions: $D_1,D_2$ and $D_3$}\label{fig1}
\end{figure}

From Lemma 1, conclutsions (b)-(d) and Theorem 1, it is easy to see that the curves $C_1:\,\alpha_2=\alpha_3$ and $C_2:\,\frac{\alpha_2}{\alpha_3}=\cos\left(\pi\big(1-\sqrt{\frac{\alpha_3^2-\alpha_2^2}{4+\alpha_3^2-\alpha_2^2}}\big)\right)$ divide the first quadrant of the $(\alpha_2,\alpha_3)$ parameter plane into three regions $D_1,D_2$ and $D_3$ (see Figure 1).

When $(\alpha_2,\alpha_3)\in D_1$, the zero solution of \eqref{22} is asymptotically stable for all $\tau\geqslant0$. When $(\alpha_2,\alpha_3)\in D_2$, the zero solution of \eqref{22} is conditionally stable with respect to the time delay $\tau$. When $(\alpha_2,\alpha_3)\in D_3\cup C_2$, the zero solution of \eqref{22} is unstable. Furthermore, when $(\alpha_2,\alpha_3)\in D_2\cup D_3\cup C_2$ (i.e. $\alpha_3>\alpha_2$), the system \eqref{22} possibly undergoes Hopf bifurcation or double Hopf bifurcation regarding the time delay $\tau$ as a bifurcation parameter. In particular,\\
{\bf Case \Rmnum{1}}: $\tau_k^+\neq\tau_l^-$ for some integers $k,l\geqslant0$ and given $(\alpha_2,\alpha_3)$ satisfying $\alpha_3>\alpha_2$, the system \eqref{22} possibly undergoes a Hopf bifurcation at $\tau=\tau_j^+$ or $\tau_j^-,j=0,1,2,\cdots$, regarding the time delay $\tau$ as the bifurcation parameter.\\
{\bf Case \Rmnum{2}}: $\tau_k^+=\tau_l^-:=\tau_0$ for some integers $k,l\geqslant0$, $\alpha_3=\alpha_{30}$ and given $\alpha_2$ satisfying $\alpha_2<\alpha_{30}$, the system \eqref{22} possibly undergoes a double Hopf bifurcation at $(\tau,\alpha_3)=(\tau_0,\alpha_{30})$ regarding the time delay $\tau$ and $\alpha_3$ as the bifurcation parameters. (Of course, we can also regard $\tau$ and $\alpha_2$ as the bifurcation parameters instead of $\tau$ and $\alpha_3$). As shown in Figure 2, there are many intersection points of $\tau_k^+$ and $\tau_l^-$ for integers $k,l\geqslant0$ and $\alpha_3$ with fixed $\alpha_2$.

\begin{figure}[H]
  \centering
  \subfigure[$\alpha_2=0.5$]{
  \includegraphics[width=3in]{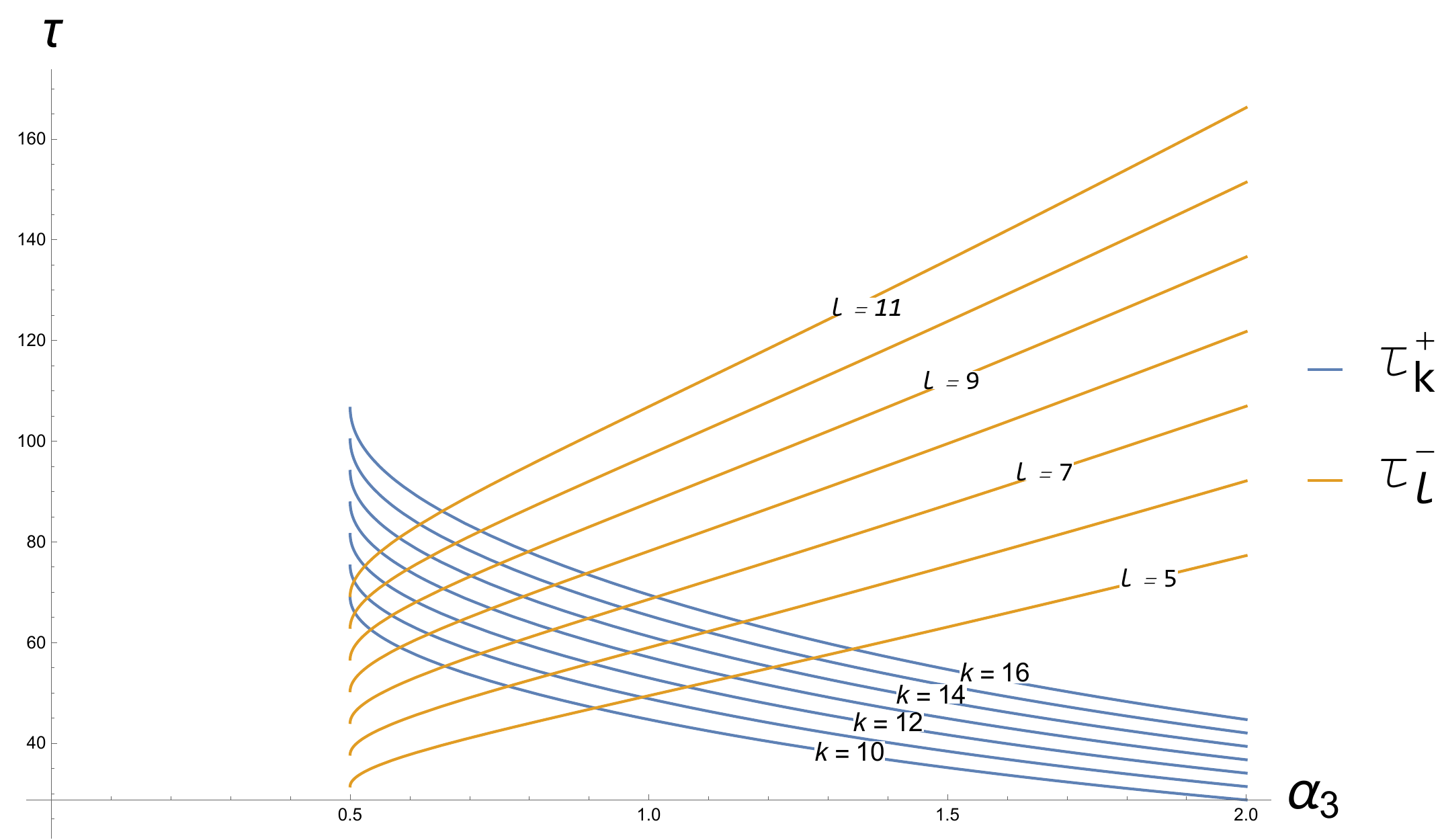}}
  \subfigure[$\alpha_2=1$]{
  \includegraphics[width=3in]{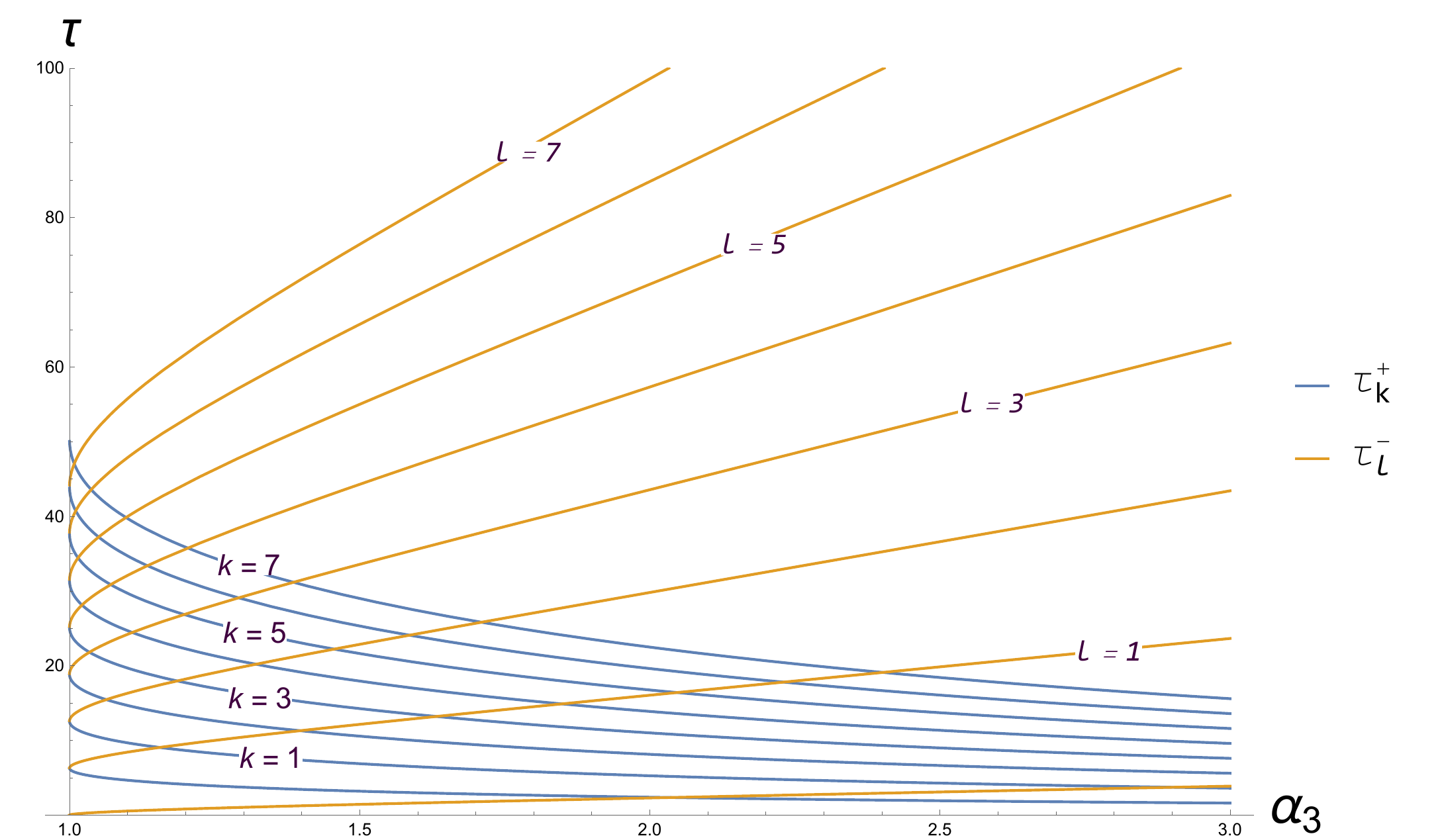}}
  \caption{The Hopf bifurcation curves $\tau=\tau_k^+(\alpha_3)$ and $\tau=\tau_l^-(\alpha_3)$ of \eqref{22} for (a) $\alpha_2=0.5$ and (b) $\alpha_2=1$. The intersection points of $\tau_k^+$ and $\tau_l^-$ are possibly double Hopf bifurcation points where the characteristic equation \eqref{23} has two pairs of simple imaginary roots $\pm{\rm i}\omega_+$ and $\pm{\rm i}\omega_-$.}\label{fig2}
\end{figure}

\noindent{\bf Case \Rmnum{3}}: $\alpha_3=\alpha_2$, because the stability of the zero solution of \eqref{22} does not change near $\tau=\tau_j(j=0,1,2,\cdots)$, the time delay $\tau$ can no longer be regarded as the Hopf bifurcation parameter. We will discuss the Hopf bifurcation preblem in \eqref{22} with $\tau=\tau_j(j=1,2,\cdots)$ regarding $\alpha_3-\alpha_2$ as a bifurcation parameter.

In the subsequent sections, we will analyze these bifurcations for {\bf Cases \Rmnum{1}-\Rmnum{3}} by using the theory of delayed differential equations.

\subsection{Decomposition of DDEs with parameters near a bifurcation value}
In this subsection, we give a decomposition of DDEs with parameters into a combination of a finite-dimensional `center' part and an infinite-dimensional `hyperbolic' part which will be used to analyze these bifurcations for {\bf Cases \Rmnum{1}-\Rmnum{3}}. Consider the following delayed differential equation
\begin{equation}\label{210}
\dot x(t)=L(\mu)x_t+F(x_t,\mu),
\end{equation}
where $x\in\mathbb R^n$, $x_t\in\mathcal C:=C\big([-1,0],\mathbb R^n\big)$ ($\mathcal C$ is the phase space of \eqref{210}, with the supremum norm), $x_t(\theta)=x(t+\theta)$ for $-1\leqslant\theta\leqslant0$, the parameter $\mu\in\mathbb R^p$ is in a neighborhood $N_p$ of the origin, $L(\mu):\mathcal C\rightarrow\mathbb R^n$ is a linear operator which is $C^k$ in $\mu$, and $F:\mathcal C\times N_p\rightarrow\mathbb R^n$ is a $C^k(k\geqslant2)$ vector function with $F(0,\mu)=0,D_1F(0,\mu)=0$ for all $\mu\in N_p$, $D_1F$ represents the derivative of $F$ with respect to the first variable $x_t$. Set
\begin{equation*}
L_0=L(0),\qquad G(x_t,\mu)=\big(L(\mu)-L_0\big)x_t+F(x_t,\mu),
\end{equation*}
then \eqref{210} can be rewritten as
\begin{equation}\label{211}
\dot x(t)=L_0x_t+G(x_t,\mu).
\end{equation}
By the Riesz representation theorem, there is a matrix $\eta(\theta)$ whose elements are bounded variation functions in $\theta\in[-1,0]$ such that
\begin{equation}\label{211+}
L_0\phi=\int_{-1}^0{\rm d}\eta(\theta)\phi(\theta),\qquad \forall\phi\in\mathcal C.
\end{equation}
We define a linear operator
\begin{equation*}
U\phi=\frac{{\rm d}\phi}{{\rm d}\theta},\qquad D(U)=\left\{\phi\in\mathcal C:\, \frac{{\rm d}\phi}{{\rm d}\theta}\in\mathcal C,\left.\frac{{\rm d}\phi}{{\rm d}\theta}\right|_{\theta=0}=L_0\phi\right\}
\end{equation*}
and a function matrix
\begin{equation*}
X_0(\theta)=\left\{\begin{aligned}
&\mathbb I_n,\quad \theta=0,\\
&0,\quad \theta\in[-1,0),
\end{aligned}\right.
\end{equation*}
where $D(U)$ is the domain of $U$, $\mathbb I_n$ is the $n\times n$ identity matrix. The spectrum of $U$ consists of all characteristic roots of the linear system $\dot x(t)=L_0x_t$. In fact, $U$ is the infinitesimal generator of this linear system.

Let $\mathcal C^\ast:=C\big([0,1],\mathbb R^{n^\ast}\big)$, where $R^{n^\ast}$ is the n-dimensional row-vector space. Define a bilinear form $\langle\cdot,\cdot\rangle$ on $\mathcal C^\ast\times\mathcal C$ by
\begin{equation*}
\langle\psi,\phi\rangle:=\psi(0)\phi(0)-\int_{-1}^0\!\!\int_{\xi=0}^\theta\psi(\xi-\theta){\rm d}\eta(\theta)\phi(\xi){\rm d}\xi,\qquad \psi\in\mathcal C^\ast,\phi\in\mathcal C
\end{equation*}
and the formal adjoint operator $U^\ast$ of $U$ by
\begin{equation*}
\langle\psi,U\phi\rangle=\langle U^\ast\psi,\phi\rangle,\qquad \phi\in D(U),\psi\in D(U^\ast).
\end{equation*}
Then we can imply
\begin{equation*}
U^\ast\psi=-\frac{{\rm d}\psi}{{\rm d}s},\qquad D(U^\ast)=\left\{\phi\in\mathcal C^\ast:\, \frac{{\rm d}\phi}{{\rm d}s}\in\mathcal C^\ast,\left.\frac{{\rm d}\phi}{{\rm d}s}\right|_{s=0}=-\int_{-1}^0\psi(-\theta){\rm d}\eta(\theta)\right\}.
\end{equation*}
Let $\land_0$ represent the set consisting of all eigenvalues of $U$ with a zero real part, $\Phi(\theta)$ be the base of the generalized eigenspace $\mathscr M_{\land_0}(U)$ corresponding to $\land_0$ and $\Psi(s)$ be the adjoint base satisfying $\langle\Psi,\Phi\rangle=\mathbb I_k$\big($k=$ the dimensional of $\mathscr M_{\land_0}(U)$\big). Then the phase space $\mathcal C$ is decomposed by $\land_0$ as $\mathcal C=P\oplus Q$, where
\begin{equation*}
P={\rm span}\Phi,\qquad Q=\{\phi\in\mathcal C:\,\langle\psi,\phi\rangle=0\,\text{ for all }\,\psi\in{\rm span}\Psi\}.
\end{equation*}
Thus, each $\phi\in\mathcal C$ can be expressed as
\begin{equation*}
\phi=\Phi\langle\Psi,\phi\rangle+\phi^Q,\qquad \phi^Q\in Q
\end{equation*}
and the solution $x_t$ of \eqref{211} can be decomposed as
\begin{equation*}
x_t(\theta)=\Phi(\theta)\langle\Psi,x_t\rangle+V(\theta,t):=\Phi(\theta)z(t)+V(\theta,t),
\end{equation*}
where $V(\cdot,t)\in Q_1=:\{\phi\in Q:\,\frac{{\rm d}\phi}{{\rm d}\theta}\in\mathcal C\}$. Hence, \eqref{211} can be decomposed as
\begin{equation}\label{212}
\left\{\begin{aligned}
&\dot z=J_0z+\Psi(0)G\big(\Phi z+V(\theta,t),\mu\big),\\
&\frac{{\rm d}V(\theta,t)}{{\rm d}t}=U_{Q_1}V(\theta,t)+(\mathbb I-\Pi)X_0(\theta)G\big(\Phi z+V(\theta,t),\mu\big),
\end{aligned}\right.
\end{equation}
where $J_0=\langle\Psi,U\Phi\rangle$ (or $J_0$ satisfies $U\Phi=\Phi J_0$), $U_{Q_1}=U|_{Q_1}$ and
\begin{equation}\label{213}
(\mathbb I-\Pi)X_0(\theta)=\left\{\begin{aligned}
&-\Phi(\theta)\Psi(0),\quad \theta\in[-1,0),\\
&\mathbb I_n-\Phi(0)\Psi(0),\quad \theta=0.
\end{aligned}\right.
\end{equation}

\section{Hopf bifurcation analysis}
In Subsection 2.1, we have obtained some critical values of $\tau$ and $\alpha_3$ at which \eqref{22} undergoes Hopf bifurcation or double Hopf bifurcation. In this section, we discuss the direction, stability of the bifurcating periodic solutions for Hopf bifurcation cases: {\bf Case \Rmnum{1}} and {\bf Case \Rmnum{3}}, listed at the end of Subsection 2.1, regarding the time delay $\tau$ and $\alpha_3-\alpha_2$ as the bifurcation parameters, respectively.

\subsection{{\bf Case \Rmnum{1}}: $\tau$ being regarded as the bifurcation parameter}
In this subsection, we discuss the Hopf bifurcation problem in \eqref{22} regarding the delay $\tau$ as the bifurcation parameter. Suppose that $\tau_0=\tau_j^+$ or $\tau_j^-$,$(j=0,1,2,\cdots)$ \big(defined by \eqref{26+}\big) is a Hopf bifurcation value, that is, for given $\alpha_2$ and $\alpha_3$ with $\alpha_3>\alpha_2$, the characteristic equation \eqref{23} has only a pair of simple imaginary roots $\pm\omega_0$ if $\tau=\tau_0$, $\omega_0=\omega_+$ for $\tau=\tau_j^+$ and $\omega_0=\omega_-$ for $\tau=\tau_j^-$ \big($\omega_\pm$ defined in \eqref{26}\big).

By rescaling the time $t\rightarrow\tau t$ to normalize the delay and phase space, \eqref{22} can be written as the form
\begin{equation}\label{30}
\left\{\begin{aligned}
&\dot x_1(t)=\tau\dot x_2(t),\\
&\dot x_2(t)=-\tau x_1(t)-\tau\alpha_2x_2(t)+\tau\alpha_3x_2(t\!-\!1)+\tau\alpha_4x_2^3(t\!-\!1)+\tau\alpha_5x_2^5(t\!-\!1)+\cdots.
\end{aligned}\right.
\end{equation}
Letting $\tau=\tau_0+\mu$, represent \eqref{30} in the form of \eqref{210}
\begin{equation}\label{31}
\dot x(t)=L(\mu)x_t+F(x_t,\mu),
\end{equation}
where $x=(x_1,x_2)^\top$,
\begin{align}
&L(\mu)x_t=(\tau_0+\mu)\big(A_1x_t(0)+A_2x_t(-1)\big),\nonumber\\
&A_1=\begin{pmatrix}0&1\\-1&-\alpha_2\end{pmatrix},\qquad A_2=\begin{pmatrix}0&0\\0&\alpha_3\end{pmatrix},\nonumber\\
&F(x_t,\mu)=(\tau_0+\mu)\dbinom{0}{\alpha_4x_2^3(t-1)+\alpha_5x_2^5(t-1)+{\rm h.o.t.}}.\label{32}
\end{align}
Based on these notations defined in Subsection 2.2 corresponding to \eqref{31}, by direct computation, we obtain
\begin{equation*}
\eta(\theta)=\left\{\begin{aligned}
&\tau_0A_1,\quad \theta=0,\\
&-\tau_0A_2\delta(\theta+1),\quad \theta\in[-1,0),
\end{aligned}\right.
\end{equation*}
where $\delta(\theta)$ is the Dirac delta function, $\land_0=\{\pm{\rm i}\tau_0\omega_0\}$ and the corresponding bases are
\begin{equation*}
\Phi(\theta)=\big(q(\theta),\bar q(\theta)\big),\qquad \Psi(s)=\dbinom{q^\ast(s)}{\overline{q^\ast}(s)}:=\big(q^\ast(s),\overline{q^\ast}(s)\big)^\top
\end{equation*}
satisfying $\langle\Psi,\Phi\rangle=\mathbb I_2$, where the notation $\bar u$ represents the complex conjugation of $u$,
\begin{align*}
&q(\theta)=\dbinom{1}{{\rm i}\omega_0}e^{{\rm i}\tau_0\omega_0\theta}, \quad \theta\in[-1,0],\\
&q^\ast(s)=d_0(1,-{\rm i}\omega_0)e^{-{\rm i}\tau_0\omega_0s}, \quad s\in[0,1],\\
&d_0=\left(1+\omega_0^2+\tau_0\alpha_3\omega_0^2e^{-{\rm i}\tau_0\omega_0}\right)^{-1}.
\end{align*}
Thus, based on the decomposition
\begin{equation}\label{33-}
x_t(\theta)=\Phi(\theta)z(t)+V(\theta,t),
\end{equation}
the system \eqref{31} is decomposed as
\begin{equation}\label{33}
\left\{\begin{aligned}
&\dot z=J_0z+\Psi(0)G\big(\Phi z+V(\theta,t),\mu\big),\\
&\frac{{\rm d}V(\theta,t)}{{\rm d}t}=U_{Q_1}V(\theta,t)+(\mathbb I-\Pi)X_0(\theta)G\big(\Phi z+V(\theta,t),\mu\big),
\end{aligned}\right.
\end{equation}
where $z=(z_1,z_2)^\top$ with $z_2=\bar z_1$, $J_0={\rm diag}({\rm i}\tau_0\omega_0,-{\rm i}\tau_0\omega_0)$ and
\begin{equation*}
G(x_t,\mu)=\big(L(\mu)-L(0)\big)x_t+F(x_t,\mu).
\end{equation*}

\vskip 0.2in

By the center manifold theorem of delayed differential equations depending on parameters (note: $F$ begins with terms of order 3), we can obtain the center manifold of \eqref{33}
\begin{equation}\label{33+}
V(\theta,t)=H(\theta,z,\mu),
\end{equation}
where
\begin{equation*}
H(\theta,z,\mu)=H_{10}(\theta,\mu)z_1+H_{01}(\theta,\mu)\bar z_1+O(|z|^3)
\end{equation*}
with $H_{10}(\theta,0)=H_{01}(\theta,0)=0$.

The system \eqref{33} restricted on the center manifold is of the following form
\begin{equation}\label{34}
\dot z=J(\mu)z+\sum_{j+k=3}\widetilde{f_{jk}}(\mu)z_1^j\bar z_1^k+O(|z|^5),
\end{equation}
where
\begin{equation*}
J(\mu)=J_0+\mu\begin{pmatrix}{\rm i}\omega_0d_0(1+\omega_0^2)&-{\rm i}\omega_0d_0(1-\omega_0^2)\\{\rm i}\omega_0\bar d_0(1-\omega_0^2)&-{\rm i}\omega_0\bar d_0(1+\omega_0^2)\end{pmatrix}+O(\mu^2),
\end{equation*}
which follows from the fact that $\pm{\rm i}\omega_0$ are roots of \eqref{23},
\begin{equation*}
\widetilde{f_{jk}}(0)=\Psi(0)F_{jk},\qquad j+k=3,
\end{equation*}
and $F_{jk}$ are the coefficients of $F$ in $z_1$ and $\bar z_1$,
\begin{equation*}
F\big(\Phi z+H(\cdot,z,0),0\big)=\sum_{j+k=3}F_{jk}z_1^j\bar z_1^k+O(|z|^5).
\end{equation*}

We only need to compute $\widetilde{f_{21}}(0)$ used later as follows
\begin{equation*}
F_{21}={\rm i}3\tau_0\alpha_4\omega_0^3e^{-{\rm i}\tau_0\omega_0}\dbinom{0}{1},\qquad \widetilde{f_{21}}(0)=3\tau_0\alpha_4\omega_0^4e^{-{\rm i}\tau_0\omega_0}\dbinom{d_0}{-\bar d_0}.
\end{equation*}
Let $\lambda_{1,2}(\mu)$ be eigenvalues of $J(\mu)$, then
\begin{align*}
\lambda_{1,2}(\mu)&=\frac{{\rm i}\mu}{2}(1+\omega_0^2)\omega_0(d_0-\bar d_0)\pm{\rm i}\tau_0\omega_0\left(1+\frac{\mu}{2\tau_0}(1+\omega_0^2)(d_0+\bar d_0)\right)+O(\mu^2)\\
&:=\sigma(\mu)\pm{\rm i}\omega(\mu).
\end{align*}
Noting the following fact
\begin{equation}\label{35}
\frac{1}{\omega_0^2}+\omega_0^2=2+\alpha_3^2-\alpha_2^2,\quad \text{i.e.}\quad (\omega_0^2-1)^2=(\alpha_3^2-\alpha_2^2)\omega_0^2,
\end{equation}
we further calculate to obtain
\begin{align*}
\sigma(\mu)&=\frac{1}{d_c}\tau_0\omega_0^2(\omega_0^2-1)(\omega_0^2+1)\mu+O(\mu^2)\\
&=\frac{\tau_0(\omega_0^2-1)(\omega_0^2+1)\mu}{4+\alpha_3^2-\alpha_2^2+2\alpha_2\tau_0+(\tau_0^2\alpha_3^2+2\alpha_2\tau_0)\omega_0^2}+O(\mu^2)\\
\omega(\mu)&=\tau_0\omega_0+\frac{1}{d_c}\omega_0(1+\omega_0^2)(1+\omega_0^2+\alpha_2\tau_0\omega_0^2)\mu+O(\mu^2)\\
&=\tau_0\omega_0+\frac{(1+\omega_0^2)(1+\omega_0^2+\alpha_2\tau_0\omega_0^2)\mu}{\omega_0(4+\alpha_3^2-\alpha_2^2+2\alpha_2\tau_0+(\tau_0^2\alpha_3^2+2\alpha_2\tau_0)\omega_0^2)}+O(\mu^2)
\end{align*}
with
\begin{equation*}
d_c=(1+\omega_0^2+\alpha_2\tau_0\omega_0^2)^2+\tau_0^2\omega_0^2(\omega_0^2-1)^2>0.
\end{equation*}
By a linear change
\begin{equation}\label{36}
z\rightarrow\big(\mathbb I_2+T(\mu)\big)z,\qquad T(0)=0
\end{equation}
(approximating to the identity transformation) to diagonalize the linear part in \eqref{34}, the first component of \eqref{34} is transformed into
\begin{equation}\label{37}
\dot z_1=\lambda_1z_1+\sum_{j+k=3}f_{1,jk}(\mu)z_1^j\bar z_1^k+O(|z|^5),
\end{equation}
where the coefficients satisfy
\begin{equation*}
f_{jk}(0)=\widetilde{f_{jk}}(0),\qquad f_{jk}(\mu)=\big(f_{1,jk}(\mu),f_{2,jk}(\mu)\big)^\top,\qquad j+k=3.
\end{equation*}
Note that we really need to study the first component, since the second component is simply the complex conjugate of the first one, and we retain the same variables $z=(z_1,\bar z_1)^\top$.

By an invertible parameter-dependent change of complex coordinate
\begin{equation*}
z_1\rightarrow z_1+a_{30}z_1^3+a_{12}z_1\bar z_1^2+a_{03}\bar z_1^3
\end{equation*}
again, \eqref{37} can be transformed into the following equation with only the resonant cubic term
\begin{equation}\label{38}
\dot z_1=\lambda_1z_1+f_{1,21}(\mu)z_1^2\bar z_1+O(|z|^5)
\end{equation}
with
\begin{equation*}
f_{1,21}(0)=\frac{3\tau_0\alpha_4}{\alpha_3d_c}\omega_0^4(\alpha_2+\omega_0^2+\tau_0\alpha_3^2\omega_0^2)+{\rm i}\frac{3\tau_0\alpha_4}{\alpha_3d_c}\omega_0^3(\omega_0^2+1)(\omega_0^2-1).
\end{equation*}
As $\alpha_2>0,\alpha_3>0$ and $\alpha_4<0$,
\begin{equation*}
{\rm Re}f_{1,21}(0)<0,\qquad \left.\frac{{\rm d}}{{\rm d}\mu}({\rm Re}\lambda_1)\right|_{\mu=0}=\sigma^\prime(0)=\frac{1}{d_c}\tau_0\omega_0^2(\omega_0^2+1)(\omega_0^2-1).
\end{equation*}
The Hopf bifurcation theory (see \cite{Wig03,Kuz98} for example) implies the following lemma.
\begin{lem}\label{lem4}
If $\tau_0=\tau_j^+$(or $\tau_j^-$), $j=0,1,2,\cdots$, then for sufficiently small $\mu>0$ (or $-\mu>0$), \eqref{38} possesses an asymptotically stable periodic solution with an amplitude $r(\mu)$ and a frequency $\omega(\mu)$ satisfying $r(\mu)=\sqrt{\frac{-\sigma(\mu)}{{\rm Re}f_{21}^\prime(0)}}+O(|\mu|^\frac{3}{2})$ and $\omega(\mu)=\tau_0\omega_0+O(\mu)$, respectively.
\end{lem}
By Theorem \ref{thm1}, Lemma \ref{lem4}, \eqref{33-}, \eqref{33+} and \eqref{36}, we obtain the following result on the Hopf bifurcation for system \eqref{22}.
\begin{thm}\label{thm2}
For fixed $(\alpha_2,\alpha_3)$ satisfying the conditions in {\bf Case \Rmnum{1}}, system \eqref{22} exhibits a Hopf bifurcation at the origin when $\tau=\tau_j^\pm$, which is subcritical for $\tau=\tau_j^-$ and supercritical for $\tau=\tau_j^+$, $j=0,1,2,\cdots$. Moreover, the bifurcating periodic solutions are asymptotically stable for the case with \eqref{27} and $j=0,1,\cdots,m$; and unstable for the case with \eqref{27} and $j\geqslant m+1$ or for the cases with \eqref{28} or \eqref{29}, in which $m$ is defined in Lemma \ref{lem3}. And bifurcating periodic solutions have following approximate expressions
\begin{equation*}
\dbinom{x_1(t)}{x_2(t)}=a(\mu)\dbinom{\cos\big(t\omega(\mu)\big)}{-\omega_0\sin\big(t\omega(\mu)\big)}+O(|\mu|^\frac{3}{2}),
\end{equation*}
where
\begin{equation*}
a(\mu)=\frac{2}{k_3\omega_0}\sqrt{\frac{(\omega_0^2+1)(\omega_0^2-1)\mu}{\alpha_2+\omega_0^2+\tau_0\alpha_3^2\omega_0^2}},\qquad \omega(\mu)=\omega_0+O(\mu),
\end{equation*}
in which $\omega_0=\omega_\pm$ and $\mu=(\tau-\tau_j^\pm)$ if $\tau_0=\tau_j^\pm$, $\tau>\tau_j^+$ for $\tau_0=\tau_j^+$ and $\tau<\tau_j^-$ for $\tau_0=\tau_j^-$, $j=0,1,2,\cdots$.
\end{thm}

\subsection{{\bf Case \Rmnum{3}}: $\alpha_3-\alpha_2$ being regarded as the bifurcation parameter}
In this subsection, we will fix $\tau=\tau_0:=2\pi j(j=1,2,\cdots)$ and regard $\alpha_3-\alpha_2$ as the Hopf bifurcation parameter with the critical value $\alpha_3-\alpha_2=0$, where the characteristic equation \eqref{23} has a pair of simple imaginary roots $\pm{\rm i}$. By rescaling the time $t\rightarrow\tau_0t$, and letting $\alpha_3=\alpha_2+\mu$, \eqref{22} can be written as
\begin{equation}\label{39}
\dot x(t)=L(\mu)x_t+F(x_t,\mu),
\end{equation}
where $x=(x_1,x_2)^\top$ and
\begin{align*}
&L(\mu)x_t=\tau_0\big(A_1(\mu)x_t(0)+A_2(\mu)x_t(-1)\big),\\
&A_1(\mu)=\begin{pmatrix}0&1\\-1&-\alpha_2\end{pmatrix},\qquad A_2(\mu)=\begin{pmatrix}0&0\\0&\alpha_2+\mu\end{pmatrix},\\
&F(x_t,\mu)=\tau_0\dbinom{0}{\alpha_4x_2^3(t-1)+\alpha_5x_2^5(t-1)+{\rm h.o.t.}}.
\end{align*}
Using these notations defined in Subsection 2.2 for system \eqref{39}, we have
\begin{equation*}
\eta(\theta)=\left\{\begin{aligned}
&\tau_0A_1(0),\quad \theta=0,\\
&-\tau_0A_2(0)\delta(\theta+1),\quad \theta\in[-1,0),
\end{aligned}\right.
\end{equation*}
where $\delta(\theta)$ is the Dirac delta function. $\land_0=\{\pm{\rm i}\}$, and corresponding bases are
\begin{equation*}
\Phi(\theta)=\big(q(\theta),\bar q(\theta)\big),\qquad \Psi(s)=\big(q^\ast(s),\overline{q\ast}(s)\big)^\top
\end{equation*}
satisfying $\langle\Psi,\Phi\rangle=\mathbb I_2$, where
\begin{align*}
&q(\theta)=\dbinom{1}{{\rm i}}e^{{\rm i}\tau_0\theta},\quad \theta\in[-1,0],\\
&q^\ast(s)=D_0(1,-{\rm i})e^{-{\rm i}\tau_0s},\quad s\in[0,1],\qquad D_0=(2+\tau_0\alpha_2)^{-1}.
\end{align*}
Hence, based on the decomposition $x_t(\theta)=\Phi(\theta)z(t)+V(\theta,t)$ with $z=(z_1,\bar z_1)^\top$, similar to the procedure in Subsection 3.1, we can obtain the center manifold of \eqref{39}
\begin{equation*}
V(\theta,t)=H(\theta,z,\mu):=H_{10}(\theta,\mu)z_1+H_{01}(\theta,\mu)\bar z_1+O(|z|^3)
\end{equation*}
with $H_{10}(\theta,0)=0,H_{01}(\theta,0)=0$, and system \eqref{39} restricted on the center manifold is
\begin{equation}\label{310}
\dot z=J(\mu)z+\sum_{j+k=3}\widetilde{f_{jk}}(\mu)z_1^j\bar z_1^k+O(|z|^5),
\end{equation}
where
\begin{align*}
&J(\mu)=\tau_0\begin{pmatrix}{\rm i}+D_0\mu&-D_0\mu\\-D_0\mu&-{\rm i}+D_0\mu\end{pmatrix}+O(\mu^2),\\
&\widetilde{f_{jk}}(0)=\Psi(0)F_{jk},\quad j+k=3
\end{align*}
and, in particular,
\begin{equation*}
F_{21}={\rm i}3\tau_0\alpha_4\dbinom{0}{1},\qquad \widetilde{f_{21}}(0)=3\tau_0\alpha_4D_0\dbinom{1}{-1}.
\end{equation*}
The eigenvalues of $J(\mu)$ are $\lambda_{1,2}=\tau_0D_0\mu\pm{\rm i}\tau_0\sqrt{1-d_1^2\mu^2}$. By a linear change $z\rightarrow\big(\mathbb I_2+\Pi(\mu)\big)z,\Pi(0)=0$ to diagonalize the linear part in \eqref{310}, we have
\begin{equation}\label{311}
\dot z_1=\lambda_1z_1+\sum_{j+k=3}f_{1,jk}(\mu)z_1^j\bar z_1^k+O(|z|^5)
\end{equation}
with coefficients satisfying
\begin{equation*}
f_{jk}(0)=\widetilde{f_{jk}}(0), \qquad f_{jk}(0)=\big(f_{1,jk}(0),f_{2,jk}(0)\big),\qquad j+k=3.
\end{equation*}
Noting $D_0>0$ and $\alpha_4<0$, it implies ${\rm Re}f_{21}^1(0)=3\tau_0\alpha_4D_0<0$ and $\frac{{\rm d}}{{\rm d}\mu}({\rm Re}\lambda_1)=\tau_0D_0>0$. Hence, by the theory of Hopf bifurcation and Theorem \ref{thm1}, we obtain the following theorem on the Hopf bifurcation for system \eqref{22}.
\begin{thm}\label{thm3}
For fixed $\tau=2\pi j(j=1,2,\cdots)$, system \eqref{22} exhibits a supercritical Hopf bifurcation at the origin when $\alpha_3-\alpha_2=0$ and the bifurcating periodic solution is asymptotically stable.
\end{thm}

\section{Double Hopf bifurcation and quasi-periodic solutions}
In this section, we study the double Hopf bifurcation in system \eqref{22} for {\bf Case \Rmnum{2}}, that is, there are integers $k,l\geqslant0$, $\alpha_3=\alpha_{30}>\alpha_2$ such that $\tau_k^+=\tau_l^-:=\tau_0$ (there are many such points $(\tau_0,\alpha_{30})$, see Figure 2), and the characteristic equation \eqref{23} with $(\tau,\alpha_3)=(\tau_0,\alpha_{30})$ has two pairs of simple imaginary roots $\pm{\rm i}\omega_{10}$ and $\pm{\rm i}\omega_{20}$, which satisfy
\begin{align*}
\omega_{10}&=\omega_+(\alpha_{30}):=\frac{1}{2}(\sqrt{4+\alpha_{30}^2-\alpha_2^2}+\sqrt{\alpha_{30}^2-\alpha_2^2}),\\
\omega_{20}&=\omega_-(\alpha_{30}):=\frac{1}{2}(\sqrt{4+\alpha_{30}^2-\alpha_2^2}-\sqrt{\alpha_{30}^2-\alpha_2^2}).
\end{align*}
We list some equalities of $\omega_{10}$ and $\omega_{20}$ in Appendix A, which will be used later. Thus, we will analyze the double Hopf bifurcation regarding $(\tau,\alpha_3)$ as the bifurcation parameters at the critical value $(\tau_0,\alpha_{30})$. Then, based on the fact that the truncated normal form system, obtained in the process of the double Hopf analysis, has quasi-periodic solutions, we will discuss the existence of quasi-periodic solutions of system \eqref{22}.

\subsection{Double Hopf bifurcation analysis}
Normalizing the delay by rescaling the time $t\rightarrow\tau t$ and letting $\tau=\tau_0+\mu_1,\alpha_3=\alpha_{30}+\mu_2$ and $\mu=(\mu_1,\mu_2)$, \eqref{22} can be written as the form of \eqref{210}, that is
\begin{equation}\label{41}
\dot x(t)=L(\mu)x_t+F(x_t,\mu),
\end{equation}
where $x=(x_1,x_2)^\top$ and
\begin{align*}
&L(\mu)x_t=(\tau_0+\mu_1)\big(A_1(\mu)x_t(0)+A_2(\mu)x_t(-1)\big),\\
&A_1(\mu)=\begin{pmatrix}0&1\\-1&-\alpha_2\end{pmatrix},\qquad A_2(\mu)=\begin{pmatrix}0&0\\0&\alpha_{30}+\mu_2\end{pmatrix},\\
&F(x_t,\mu)=(\tau_0+\mu_1)\dbinom{0}{\alpha_4x_2^3(t-1)+\alpha_5x_2^5(t-1)+{\rm h.o.t.}}.
\end{align*}
Using these notations defined in Subsection 2.2 for system \eqref{41}, we have
\begin{equation*}
\eta(\theta)=\left\{\begin{aligned}
&\tau_0A_1(0),\quad \theta=0,\\
&-\tau_0A_2(0)\delta(\theta+1),\quad \theta\in[-1,0),
\end{aligned}\right.\qquad \land_0=\{\pm{\rm i}\tau_0\omega_{10},\pm{\rm i}\tau_0\omega_{20}\}.
\end{equation*}
The eigenfunction of $U$ corresponding to the eigenvalue ${\rm i}\tau_0\omega_{j0}$ is
\begin{equation*}
q_j(\theta)=\dbinom{1}{{\rm i}\omega_{j0}}e^{{\rm i}\tau_0\omega_{j0}\theta},\quad \theta\in[-1,0],\quad j=1,2,
\end{equation*}
while the eigenfunction of $U^\ast$ is
\begin{equation*}
p_j(s)=d_j(1,-{\rm i}\omega_{j0})e^{-{\rm i}\tau_0\omega_{j0}s},\quad s\in[0,1],\quad j=1,2,
\end{equation*}
where $d_j=(1+\omega_{j0}^2+\tau_0\alpha_{30}\omega_{j0}^2e^{-{\rm i}\tau_0\omega_{j0}})^{-1}$. Then,
\begin{equation*}
\Phi(\theta)=\big(q_1(\theta),\bar q_1(\theta),q_2(\theta),\bar q_2(\theta)\big)\,\text{ and }\, \Psi(s)=\big(p_1(s),\bar p_1(s),p_2(s),\bar p_2(s)\big)^\top
\end{equation*}
satisfy $\langle\Psi,\Phi\rangle=\mathbb I_4$. Based on the decomposition
\begin{equation*}
x_t(\theta)=\Phi(\theta)z+V(\theta,t),
\end{equation*}
the system \eqref{41} is decomposed as
\begin{equation}\label{42}
\left\{\begin{aligned}
&\dot z=J_0z+\Psi(0)G\big(\Phi z+V(\theta,t),\mu\big),\\
&\frac{{\rm d}V(\theta,t)}{{\rm d}t}=U_{Q_1}V(\theta,t)+(\mathbb I-\Pi)X_0(\theta)G\big(\Phi z+V(\theta,t),\mu\big),
\end{aligned}\right.
\end{equation}
where $z=(z_1,\bar z_1,z_2,\bar z_2)^\top$, $J_0={\rm diag}({\rm i}\tau_0\omega_{10},-{\rm i}\tau_0\omega_{10},{\rm i}\tau_0\omega_{20},-{\rm i}\tau_0\omega_{20})$ and
\begin{equation*}
G(x_t,\mu)=\big(L(\mu)-L(0)\big)x_t+F(x_t,\mu).
\end{equation*}
By the center manifold theorem, we can suppose the center manifold of \eqref{42} is
\begin{equation*}
V(\theta,t)=H(\theta,z,\mu):=\sum_{1\leqslant|k|\leqslant5}H_k(\theta,\mu)z^k+O(|z|^7),
\end{equation*}
where $z^k=z_1^{k_1}\bar z_1^{k_2}z_2^{k_3}\bar z_2^{k_4},|k|=k_1+k_2+k_3+k_4$ and $H_k\in Q_1$. Noting that $F$ is an odd function in the space variable and begins with terms of order 3, it implies
\begin{align*}
&H_k(\theta,\mu)\equiv0\quad\text{for}\quad|k|=2,4\\
\text{and}\quad&H_k(\theta,0)=0\quad\text{for}\quad|k|=1.
\end{align*}
Moreover, it is enough that, we only need to obtain the coefficients of terms of order 3 in the center manifold when $\mu=0$, so that we can analyze double Hopf bifurcation in system \eqref{42}.

\vskip 0.2in

Let
\begin{equation}\label{43}
F(\Phi z,0)=\sum_{|k|=3,5}F_kz^k+O(|z|^7),
\end{equation}
in which the coefficients $F_k$ with $|k|=3$ are listed in Appendix B. Equation \eqref{42} implies
\begin{equation}\label{44}
U_{Q_1}H_k(\theta,0)={\rm i}\tau_0\langle k,\tilde\omega_0\rangle H_k(\theta,0)-(\mathbb I-\Pi)X_0(\theta)F_k,
\end{equation}
where $\tilde\omega_0=(\omega_{10},-\omega_{10},\omega_{20},-\omega_{20})^\top,k=(k_1,k_2,k_3,k_4)^\top$. By the definition of $U_{Q_1}$, the continuity of $\frac{{\rm d}}{{\rm d}\theta}H_k$, $H_k\in Q_1$ and \eqref{213}, \eqref{44} is converted to the following boundary value problem
\begin{equation}\label{45}
\left\{\begin{aligned}
&\frac{{\rm d}}{{\rm d}\theta}H_k(\theta,0)={\rm i}\tau_0\langle k,\tilde\omega_0\rangle H_k(\theta,0)+\Phi(\theta)\Psi(0)F_k,\\
&\frac{{\rm d}}{{\rm d}\theta}H_k(0,0)=L(0)H_k(\theta,0)+F_k,\quad\langle\Psi,H_k\rangle=0,\qquad |k|=3,5.
\end{aligned}\right.
\end{equation}
We do not solve \eqref{45} here, since the coefficients $H_k$ with $|k|=3,5$ of the center manifold will only affect the terms of order$\geqslant5$ in the restricted system on the center manifold, and our subsequent discussion do not require these coefficients.

The system \eqref{42} restricted on the center manifold is written as
\begin{equation}\label{46}
\dot z=\tilde J(\mu)z+\sum_{|k|=3,5}\widetilde{f_k}(\mu)z^k+O(|z|^7),
\end{equation}
where
\begin{equation*}
\tilde J(\mu)=J_0+J_1(\mu),\qquad \widetilde{f_k}(0)=\Psi(0)F_k\,\text{ for }\,|k|=3,5.
\end{equation*}
The expression of $J_1(\mu)$ is placed in Appendix C. Let $\lambda_1(\mu):=\sigma_1(\mu)+{\rm i}\omega_1(\mu),\bar \lambda_1,\lambda_2(\mu):=\sigma_2(\mu)+{\rm i}\omega_2(\mu)$ and $\bar\lambda_2$ be the four eigenvalues of $\tilde J(\mu)$. Then we have
\begin{align}
&\det\left.\left(\frac{\partial(\sigma_1,\sigma_2)}{\partial(\mu_1,\mu_2)}\right)\right|_{\mu=0}=\frac{\tau_0^2(\omega_{10}^2-\omega_{20}^2)^3}{\alpha_{30}\beta_1\beta_2(3\omega_{10}^4+3\omega_{20}^4+10)}\Delta_1>0,\label{46*}\\
&\det\left.\left(\frac{\partial(\omega_1,\omega_2)}{\partial(\mu_1,\mu_2)}\right)\right|_{\mu=0}=\frac{\tau_0(\omega_{10}^2-\omega_{20}^2)^3}{\alpha_{30}\beta_1\beta_2(3\omega_{10}^4+3\omega_{20}^4+10)}(2+\tau_0\alpha_2)(\omega_{10}^2+\omega_{20}^2+2)>0,\nonumber
\end{align}
where
\begin{align}
&\Delta_1=4\alpha_2+2\tau_0\alpha_{30}+\alpha_2(\alpha_{30}^2-\alpha_2^2),\nonumber\\
&\beta_j=(1+\omega_{j0}^2+\tau_0\alpha_2\omega_{j0}^2)^2+\tau_0^2\omega_{j0}^2(\omega_{j0}^2-1)^2,\quad j=1,2.\label{46+}
\end{align}
The approximate expressions of $\sigma_j(\mu)$ and $\omega_j(\mu)(j=1,2)$ are placed in Appendix D. By a linear change $z\rightarrow z+T(\mu)z,T(0)=0$ to diagonalize the linear part in \eqref{46}, we obtain
\begin{equation}\label{47}
\dot z=Jz+\sum_{|k|=3,5}f_k(\mu)z^k+O(|z|^7),
\end{equation}
where
\begin{equation*}
J={\rm diag}(\lambda_1,\bar\lambda_1,\lambda_2,\bar\lambda_2),\qquad f_k(0)=\widetilde{f_k}(0)\,\text{ for }\,|k|=3,5.
\end{equation*}

\vskip 0.2in

Assume that $\omega_{10}$ and $\omega_{20}$ are nonresonant up to order 6, i.e., $n_1\omega_{10}\neq n_2\omega_{20}$ for all integers $n_1$ and $n_2$ satisfying $1\leqslant n_1+n_2\leqslant6$, which is equivalent to
\begin{equation}\label{48}
\alpha_{30}^2-\alpha_2^2\neq\frac{(n_1-n_2)^2}{n_1n_2},\quad \text{for all }\,0<n_1<n_2\,\text{ and }\,n_1+n_2\leqslant6,\quad n_1,n_2\in\mathbb Z_+.
\end{equation}
Under the assumption \eqref{48}, it implies by the normal form method that, there exists an invertible parameter-dependent change of complex coordinate variables
\begin{equation*}
z\rightarrow z+\sum_{|k|=3,5}g_k(\mu)z^k,
\end{equation*}
which reduce system \eqref{47}, for sufficiently small $|\mu|$, into the following form
\begin{equation}\label{49}
\left\{\begin{aligned}
\dot z_1=\lambda_1(\mu)z_1+p_{11}(\mu)z_1^2\bar z_1+p_{12}&(\mu)z_1z_2\bar z_2+P_{11}(\mu)z_1^3\bar z_1^2\\
&+P_{12}(\mu)z_1^2\bar z_1z_2\bar z_2+P_{13}(\mu)z_1z_2^2\bar z_2^2+O(|z|^7),\\
\dot z_2=\lambda_2(\mu)z_2+p_{21}(\mu)z_2^2\bar z_2+p_{22}&(\mu)z_1\bar z_1z_2+P_{21}(\mu)z_1^2\bar z_1^2z_2\\
&+P_{22}(\mu)z_1\bar z_1z_2^2\bar z_2+P_{23}(\mu)z_2^3\bar z_2^2+O(|z|^7),
\end{aligned}\right.
\end{equation}
where the coefficients $p_{ij}$ satisfy
\begin{align*}
&p_{11}(0)=-{\rm i}\omega_{10}d_1F_{2,2100},\qquad p_{12}(0)=-{\rm i}\omega_{10}d_1F_{2,1011},\\
&p_{21}(0)=-{\rm i}\omega_{20}d_2F_{2,1110},\qquad p_{22}(0)=-{\rm i}\omega_{20}d_2F_{2,0021},
\end{align*}
in which $F_{2,k_1k_2k_3k_4}$ represents the second component of $F_{k_1k_2k_3k_4}$, and the calculation of $P_{ij}$ may be found in bifurcation textbooks and papers, see e.g. \cite{Kuz98}. We do not need the value of $P_{ij}$ as the double Hopf bifurcation diagrams of our system are the `simple' case.

In the polar coordinates $(r,\theta)$, i.e. $z_1=r_1e^{{\rm i}\theta_1},z_1=r_2e^{{\rm i}\theta_2}$, system \eqref{49} becomes
\begin{equation}\label{410}
\left\{\begin{aligned}
&\dot r_1=r_1\big(\sigma_1(\mu)+a_{11}(\mu)r_1^2+a_{12}(\mu)r_2^2+A_{11}r_1^4+A_{12}r_1^2r_2^2+A_{13}r_2^4\big)+O(|r|^7),\\
&\dot r_2=r_2\big(\sigma_2(\mu)+a_{21}(\mu)r_1^2+a_{22}(\mu)r_2^2+A_{21}r_1^4+A_{22}r_1^2r_2^2+A_{23}r_2^4\big)+O(|r|^7),\\
&\dot\theta_1=\omega_1(\mu)+q_{11}(\mu)r_1^2+q_{12}r_2^2+Q_{11}r_1^4+Q_{12}r_1^2r_2^2+Q_{13}r_2^4+r_1^{-1}O(|r|^7),\\
&\dot\theta_2=\omega_2(\mu)+q_{21}(\mu)r_1^2+q_{22}r_2^2+Q_{21}r_1^4+Q_{22}r_1^2r_2^2+Q_{23}r_2^4+r_2^{-1}O(|r|^7),
\end{aligned}\right.
\end{equation}
where
\begin{align*}
a_{jl}(\mu)={\rm Re}\,p_{jl}(\mu),\quad q_{jl}(\mu)&={\rm Im}\,p_{jl}(\mu),\quad j,l=1,2,\\
A_{jl}(\mu)={\rm Re}P_{jl}(\mu),\quad Q_{jl}(\mu)&={\rm Im}P_{jl}(\mu),\quad j=1,2,l=1,2,3.
\end{align*}
The values of $a_{jl}$ and $q_{jl}$ at the bifurcation point $\mu=0$ are listed in Appendix E. In particular, we have
\begin{equation}\label{411}
\frac{a_{11}(0)}{a_{12}(0)}=\frac{\omega_{10}^2}{2\omega_{20}^2},\quad \frac{a_{21}(0)}{a_{22}(0)}=\frac{2\omega_{10}^2}{\omega_{20}^2},\quad \frac{q_{11}(0)}{q_{12}(0)}=\frac{\omega_{10}^2}{2\omega_{20}^2},\quad \frac{q_{21}(0)}{q_{22}(0)}=\frac{2\omega_{10}^2}{\omega_{20}^2}.
\end{equation}
As $\alpha_4<0$ and $a_{11}(0)<0,a_{22}(0)<0$, the double Hopf bifurcation is the `simple' case (see Page 357 in \cite{Kuz98}). The topology of the bifurcation diagram to the truncated amplitude system of \eqref{410}
\begin{equation}\label{412}
\left\{\begin{aligned}
&\dot r_1=r_1\big(\sigma_1(\mu)+a_{11}(\mu)r_1^2+a_{12}(\mu)r_2^2+A_{11}r_1^4+A_{12}r_1^2r_2^2+A_{13}r_2^4\big),\\
&\dot r_2=r_2\big(\sigma_2(\mu)+a_{21}(\mu)r_1^2+a_{22}(\mu)r_2^2+A_{21}r_1^4+A_{22}r_1^2r_2^2+A_{23}r_2^4\big)
\end{aligned}\right.
\end{equation}
is independent of the terms of order 5. In particular, by Lemma 8.17 in \cite{Kuz98}, system \eqref{412} is locally topologically equivalent near the origin to the system
\begin{equation}\label{413}
\left\{\begin{aligned}
&\dot \rho_1=2\rho_1\big(\sigma_1(\mu)-\rho_1-\delta_1(\mu)\rho_2\big),\\
&\dot \rho_2=2\rho_2\big(\sigma_2(\mu)-\delta_2(\mu)\rho_1-\rho_2\big),
\end{aligned}\right.
\end{equation}
where we introduce the new phase variables according to
\begin{equation}\label{414}
\rho_1=-a_{11}r_1^2,\qquad \rho_2=-a_{22}r_2^2
\end{equation}
and
\begin{equation}\label{414+}
\delta_1(\mu)=\frac{a_{12}(\mu)}{a_{22}(\mu)},\qquad \delta_2(\mu)=\frac{a_{21}(\mu)}{a_{22}(\mu)}.
\end{equation}
Moreover, by \eqref{411} we obtain
\begin{equation}\label{415}
\delta_1(0)>0,\quad \delta_2(0)>0,\quad \delta_1(0)\delta_2(0)=4,\quad \delta_1(0)-\delta_2(0)>0.
\end{equation}
Here, we used these expressions of $a_{11}(0)$ and $a_{22}(0)$, and some equalities given in Appendix A.

\vskip 0.2in

The \eqref{46*} implies that the map $(\mu_1,\mu_2)\mapsto(\sigma_1,\sigma_2)$ is regular at $\mu=0$. Thus, we can regard $(\sigma_1,\sigma_2)$ as the bifurcation parameters of system \eqref{413}.

\begin{figure}[H]
  \centering
  \includegraphics[width=3.6in]{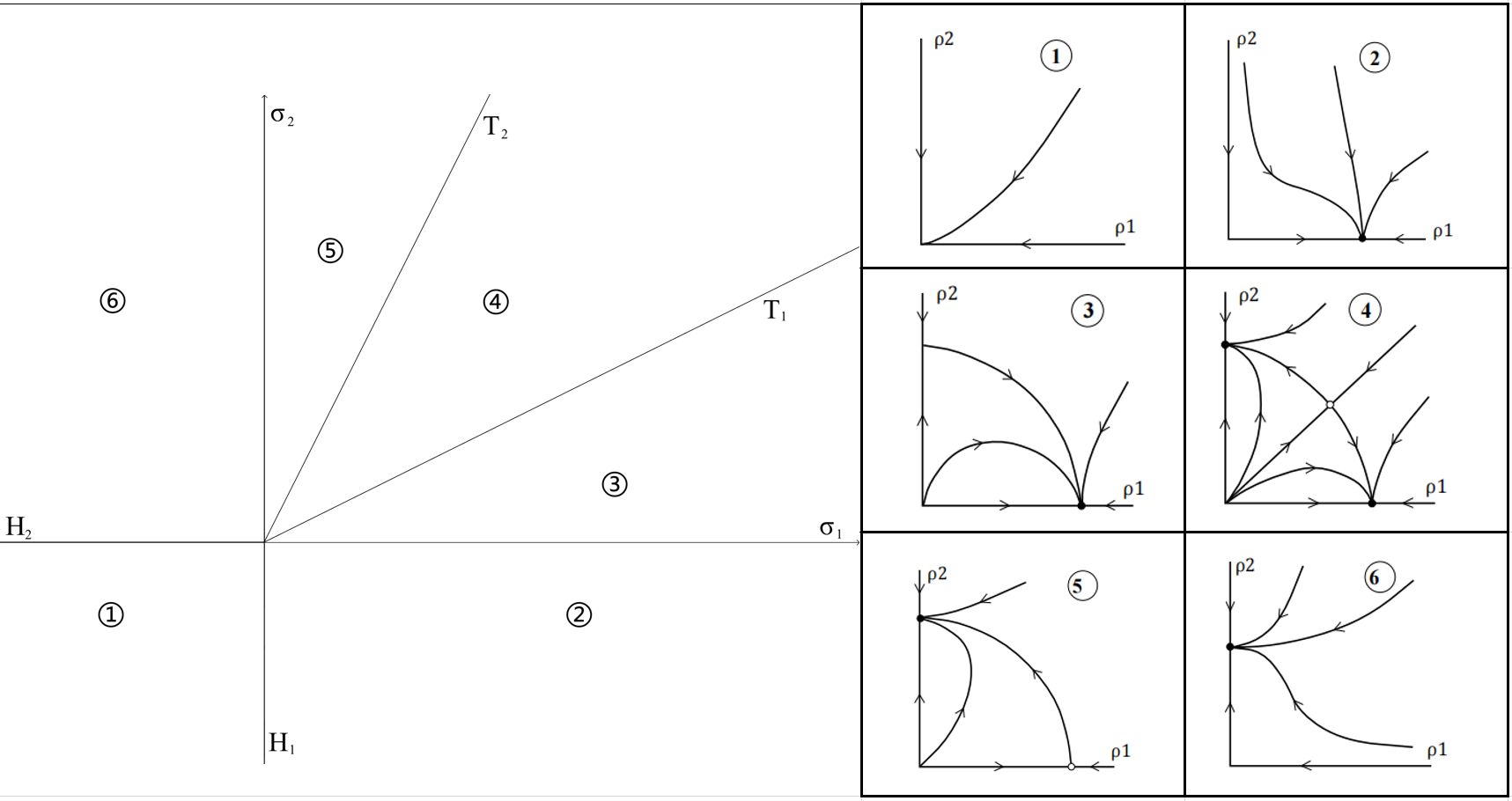}
  \caption{Bifurcation diagram of \eqref{413} regarding $\sigma_1$ and $\sigma_2$ as bifurcation parameters.}\label{fig3}
\end{figure}

System \eqref{413} has an equilibrium $E_0=(0,0)$, the origin, for all $\sigma_1$ and $\sigma_2$; two equilibria $E_1=(\sigma_1,0)$ and $E_2=(0,\sigma_2)$, bifurcating from the equilibrium $E_0$ at the bifurcation lines $H_1=\{(\sigma_1,\sigma_2):\sigma_1=0\}$ and $H_2=\{(\sigma_1,\sigma_2):\sigma_2=0\}$, respectively; an nontrivial equilibrium
\begin{equation*}
E_3=\left(-\frac{\sigma_1-\delta_1\sigma_2}{\delta_1\delta_2-1},\frac{\delta_2\sigma_1-\sigma_2}{\delta_1\delta_2-1}\right)
\end{equation*}
bifurcating from $E_1$ (or $E_2$), colliding with $E_2$ (or $E_1$) and disappearing from the positive quadrant on the bifurcation curves
\begin{equation*}
T_1=\{(\sigma_1,\sigma_2):\sigma_1=\delta_1\sigma_2,\sigma_2>0\}\quad\text{and}\quad T_2=\{(\sigma_1,\sigma_2):\sigma_2=\delta_2\sigma_1,\sigma_1>0\}.
\end{equation*}
The bifurcation diagram of system \eqref{413} is shown in Figure 3 (it is the Case \Rmnum{1} in Figure 8.25 on page 359 in \cite{Kuz98}).

Thus, for the truncated system of \eqref{410}
\begin{equation}\label{416}
\left\{\begin{aligned}
&\dot r_1=r_1\big(\sigma_1(\mu)+a_{11}(\mu)r_1^2+a_{12}(\mu)r_2^2+A_{11}r_1^4+A_{12}r_1^2r_2^2+A_{13}r_2^4\big),\\
&\dot r_2=r_2\big(\sigma_2(\mu)+a_{21}(\mu)r_1^2+a_{22}(\mu)r_2^2+A_{21}r_1^4+A_{22}r_1^2r_2^2+A_{23}r_2^4\big),\\
&\dot\theta_1=\omega_1(\mu)+q_{11}(\mu)r_1^2+q_{12}r_2^2+Q_{11}r_1^4+Q_{12}r_1^2r_2^2+Q_{13}r_2^4,\\
&\dot\theta_2=\omega_2(\mu)+q_{21}(\mu)r_1^2+q_{22}r_2^2+Q_{21}r_1^4+Q_{22}r_1^2r_2^2+Q_{23}r_2^4,
\end{aligned}\right.
\end{equation}
by the bifurcation diagram of \eqref{413} we can obtain the corresponding bifurcation of \eqref{416}: curves $H_1,H_2,T_1^\prime=\{(\sigma_1,\sigma_2):\sigma_1=\delta_1\sigma_2+O(\sigma_2^2),\sigma_2>0\}$ and $T_2^\prime=\{(\sigma_1,\sigma_2):\sigma_2=\delta_2\sigma_1+O(\sigma_1^2),\sigma_1>0\}$ divide the sufficiently small neighborhood of $(\sigma_1,\sigma_2)=(0,0)$ into six regions on the $(\sigma_1,\sigma_2)$-parameter plane depicted in Figure 4, moving anticlockwise around the origin, from region \Rmnum{1},\\

\begin{figure}[H]
  \centering
  \includegraphics[width=2in]{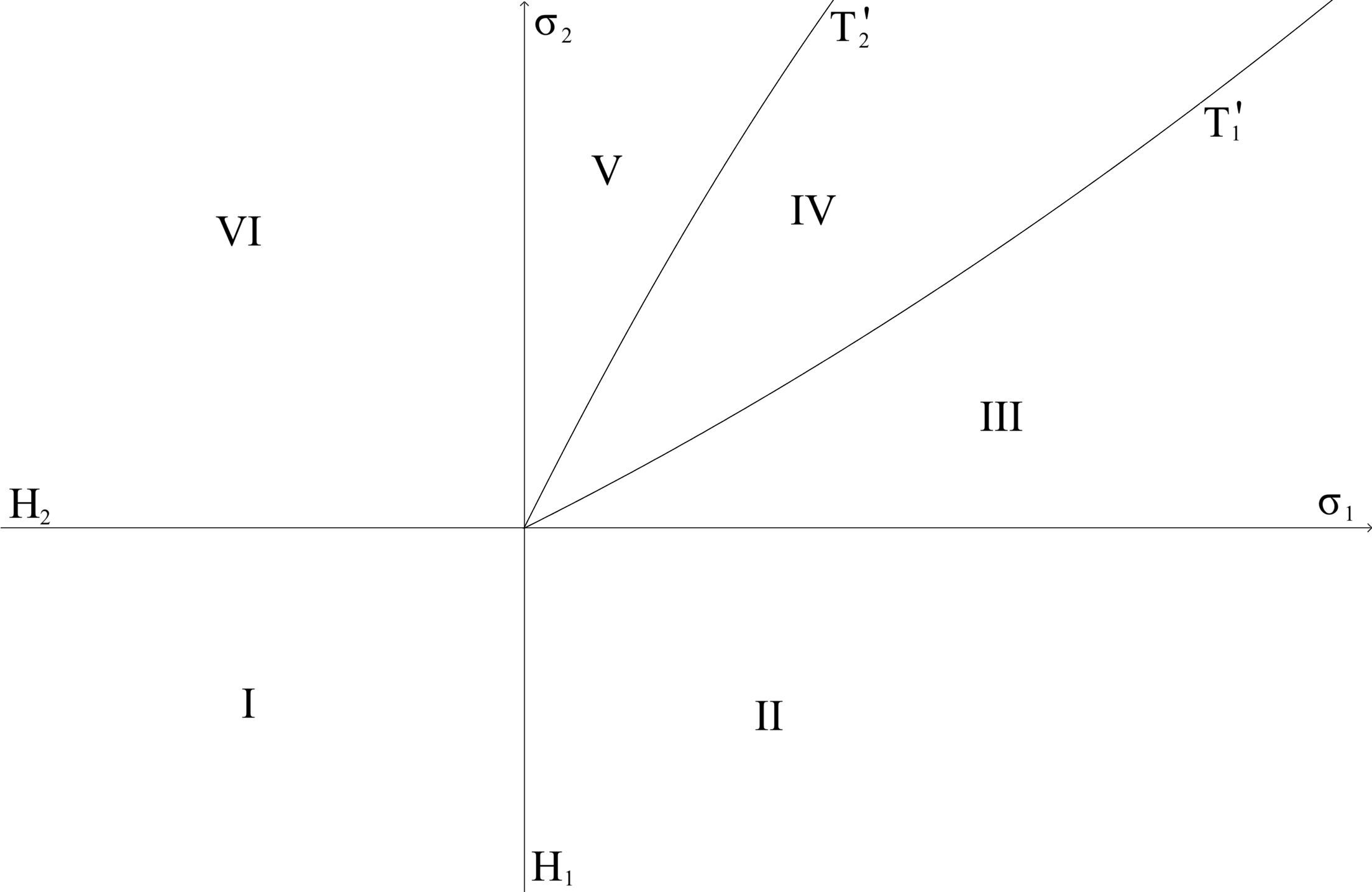}
  \caption{The parametric bifurcation portrait of \eqref{416}.}\label{fig4}
\end{figure}

\hangafter 1
\hangindent 1.5em
\noindent(\rmnum{1}) crossing the line $H_1$ from region \Rmnum{1} to region \Rmnum{2}, a stable periodic solution $\mathscr P_1$ appears by Hopf bifurcation from the trivial solution. And this stable periodic solution $\mathscr P_1$ is persistent for $(\sigma_1,\sigma_2)$ in regions \Rmnum{2}, \Rmnum{3} and \Rmnum{4};

\hangafter 1
\hangindent 1.5em
\noindent(\rmnum{2}) crossing the line $H_2$ from region \Rmnum{2} to region \Rmnum{3}, another unstable periodic solution $\mathscr P_2$ appears by Hopf bifurcation from the trivial solution, which is persistent for $(\sigma_1,\sigma_2)$ in region \Rmnum{3};

\hangafter 1
\hangindent 1.5em
\noindent(\rmnum{3}) crossing the curve $T_1^\prime$ from region \Rmnum{3} to region \Rmnum{4}, an unstable and invariant two-dimensional torus (2-torus) $\mathscr T$ appears by Neimark-Sacker bifurcation from the periodic solution $\mathscr P_2$, which is persistent for $(\sigma_1,\sigma_2)$ in regions \Rmnum{4}. Meanwhile, the periodic solution $\mathscr P_2$ becomes stable and is persistent for $(\sigma_1,\sigma_2)$ in regions \Rmnum{4}, \Rmnum{5} and \Rmnum{6};

\hangafter 1
\hangindent 1.5em
\noindent(\rmnum{4}) crossing the curve $T_2^\prime$ from region \Rmnum{4} to region \Rmnum{5}, the unstable invariant torus $\mathscr T$ collides with the stable periodic solution $\mathscr P_1$ and then disappears. Meanwhile, the periodic solution $\mathscr P_1$ becomes unstable and persistent for $(\sigma_1,\sigma_2)$ in region \Rmnum{5};

\hangafter 1
\hangindent 1.5em
\noindent(\rmnum{5}) crossing the line $H_1$ from region \Rmnum{5} to region \Rmnum{6}, the periodic solution $\mathscr P_1$ collides with the trivial solution and then disappears;

\hangafter 1
\hangindent 1.5em
\noindent(\rmnum{6}) crossing the line $H_2$ from region \Rmnum{6} to region \Rmnum{1}, the periodic solution $\mathscr P_2$ collides with the trivial solution and then disappears.

\vskip 0.2in

In particular, we have the following lemma.
\begin{lem}\label{lem5}
For each $(\sigma_1,\sigma_2)$ in region \Rmnum{4} depicted in Figure 4, system \eqref{416} possesses a hyperbolic invariant two-dimensional torus, which is either filled with quasi-periodic solutions (called quasi-periodic invariant torus) or filled with periodic solutions (called periodic invariant torus).
\end{lem}
Now, we relate the above results to the system \eqref{22}. By Theorem \ref{thm1} and its following analysis, it is easy to see that the center manifold $V(\theta,t)=H(\theta,z,0)$ of \eqref{22} is stable if $(\alpha_2,\alpha_{30})\in D_2$ in Figure 1 and $\tau_{m+1}^+=\tau_{m+1}^-(:=\tau_0)$ (i.e. $k=l=m+1$ in {\bf Case \Rmnum{2}}); is unstable if $(\alpha_2,\alpha_{30})\in D_2$ and $\tau_k^+=\tau_l^-(:=\tau_0)$ with some $k\neq l$ (then $k>l\geqslant m+1$) or $(\alpha_2,\alpha_{30})\in D_3\cup C_2$ in Figure 1 in {\bf Case \Rmnum{2}}. Thus, we obtain the following theorem.
\begin{thm}\label{thm4}
For fixed $\alpha_2$ and $(\tau_0,\alpha_{30})$ satisfying the conditions in {\bf Case \Rmnum{2}} and \eqref{48}, system \eqref{22} exhibits a `simple' double bifurcation at the origin when $(\tau,\alpha_3)=(\tau_0,\alpha_{30})$. On the $(\sigma_1,\sigma_2)$-parameter plane based on the regular map $(\mu_1,\mu_2)\mapsto(\sigma_1,\sigma_2)$, for sufficiently small $|\mu|$, there exist bifurcation curves $\tilde H_j,\tilde T_j(j=1,2)$ being tangent to $H_j,T_j^\prime(j=1,2)$ depicted in Figure 4, respectively, and at the curves $\tilde H_j,\tilde T_j(j=1,2)$, Hopf bifurcations of the origin and Neimark-Sacker bifurcations of limit cycles occur, respectively. Moreover, the bifurcating invariant two-dimensional tori at the bifurcation curves $\tilde T_j,j=1,2$, are unstable; the bifurcating periodic solutions at Hopf bifurcation curves $\tilde H_j,j=1,2$ are of the same stability of ones in system \eqref{416} only for the case where $\tau_{m+1}^+=\tau_{m+1}^-(:=\tau_0)$ and $(\alpha_2,\alpha_{30})\in D_2$, and are unstable in other cases.
\end{thm}

\subsection{The existence of quasi-periodic solutions}
Lemma \ref{lem5} shows that the truncated system \eqref{416} of system \eqref{410} possesses an invariant 2-torus for each $(\sigma_1,\sigma_2)$ in region \Rmnum{4}. Kuznetsov pointed out that, for the full system \eqref{410} (or equivalently system \eqref{22}), an invariant 2-torus only appears near Neimark-Sacker bifurcation curves $\tilde H_j,j=1,2$, and, away from these curves, the 2-torus loses smootheness and is destroyed. Moreover, the orbit structure on a torus of the full system \eqref{410} is generically different from that of the truncated system \eqref{416}. An invariant periodic torus of \eqref{416} is ruptured generically under perturbations (adding the removed higher-order terms). We will prove that, for most (in the sense of Lebesgue measure) of the parameter values in region \Rmnum{4}, an invariant quasi-periodic 2-torus remains under perturbations by using KAM theorems in \cite{LS19}. To this end, we need to reduce system \eqref{410} to a standard form suitable for applying these KAM theorems by some transformations.

Inserting the change \eqref{414} in \eqref{410}, we have
\begin{equation}\label{417}
\left\{\begin{aligned}
&\dot \rho_1=2\rho_1\big(\sigma_1-\rho_1-\delta_1\rho_2+B_{11}\rho_1^2+B_{12}\rho_1\rho_2+B_{13}\rho_2^2\big)+O(|\rho|^4),\\
&\dot \rho_2=2\rho_2\big(\sigma_2-\delta_2\rho_1-\rho_2+B_{21}\rho_1^2+B_{22}\rho_1\rho_2+B_{23}\rho_2^2\big)+O(|\rho|^4),\\
&\dot\theta_1=\tilde\omega_1(\sigma)+s_{11}\rho_1+s_{12}\rho_2+S_{11}\rho_1^2+S_{12}\rho_1\rho_2+S_{13}\rho_2^2+\rho_1^{-\frac{1}{2}}O(|\rho|^\frac{7}{2}),\\
&\dot\theta_2=\tilde\omega_2(\sigma)+s_{21}\rho_1+s_{22}\rho_2+S_{21}\rho_1^2+S_{22}\rho_1\rho_2+S_{23}\rho_2^2+\rho_2^{-\frac{1}{2}}O(|\rho|^\frac{7}{2}),
\end{aligned}\right.
\end{equation}
where $\rho=(\rho_1,\rho_2)^\top$ and
\begin{align*}
&\tilde\omega_j(\sigma)=\omega_j\big(\mu(\sigma)\big),\quad B_{j1}=\frac{A_{j1}}{a_{11}^2},\quad B_{j2}=\frac{A_{j2}}{a_{11}a_{22}},\quad B_{j3}=\frac{A_{j3}}{a_{22}^2},\\
&s_{j1}=-\frac{q_{j1}}{a_{11}},\quad s_{j2}=-\frac{q_{j2}}{a_{22}},\quad S_{j1}=\frac{Q_{j1}}{a_{11}^2},\quad S_{j2}=\frac{Q_{j2}}{a_{11}a_{22}},\quad S_{j3}=\frac{Q_{j3}}{a_{22}^2},\qquad j=1,2.
\end{align*}
Here, we regarded $\mu$ as a function in $\sigma$ by \eqref{46*}, all coefficients depend on the new parameter $\sigma$. Taking the rescaling and translational changes
\begin{equation*}
\rho_j\rightarrow\varepsilon\rho_j,\quad \sigma_j\rightarrow\varepsilon\sigma_j,\quad \rho_j=\rho_{j0}+\varepsilon^\frac{1}{2}u_j,\qquad j=1,2
\end{equation*}
in \eqref{417}, where $\varepsilon>0$ is a sufficiently small parameter and
\begin{equation}\label{417+}
\rho_{10}=-\frac{\sigma_1-\delta_1\sigma_2}{\delta_1\delta_2-1},\quad \rho_{20}=\frac{\delta_2\sigma_1-\sigma_2}{\delta_1\delta_2-1},
\end{equation}
we have
\begin{equation}\label{418}
\left\{\begin{aligned}
&\dot u=\varepsilon\big(\Omega^0(\sigma)u+\varepsilon^\frac{1}{2}f_1(u,\theta,\sigma,\varepsilon)\big),\\
&\dot \theta=\omega^0(\delta,\varepsilon)+\varepsilon^\frac{3}{2}f_2(u,\theta,\sigma,\varepsilon),
\end{aligned}\right.
\end{equation}
where $u=(u_1,u_2)^\top,\theta=(\theta_1,\theta_2)^\top$,
\begin{align}
&\Omega^0(\sigma)=\begin{pmatrix}-2\rho_{10}&-2\delta_1\rho_{10}\\-2\delta_2\rho_{20}&-2\rho_{20}\end{pmatrix},\nonumber\\
&\omega^0(\delta,\varepsilon)=\dbinom{\tilde\omega_1(\varepsilon\sigma)+\varepsilon(s_{11}\rho_{10}+s_{12}\rho_{20}}{\tilde\omega_2(\varepsilon\sigma)+\varepsilon(s_{21}\rho_{10}+s_{22}\rho_{20}}:=\dbinom{\omega_{10}}{\omega_{20}}+\varepsilon\omega^1(\sigma)+O(\varepsilon^2),\label{419}
\end{align}
and $f_1,f_2$ are sufficiently smooth functions of all variables. Let $\tilde\Pi_0$ be a closed and convex subset of region \Rmnum{4}, whose distance is at least $d_0\varepsilon$ to the boundary of region \Rmnum{4} ($d_0>0$ is a positive number such that $\tilde\Pi_0$ is of positive Lebesgue measure). Set $\Pi_0=\frac{1}{\varepsilon}\tilde\Pi_0$.
\begin{thm}\label{thm5}
Under these conditions of Theorem \ref{thm4}, for any given $0<\gamma\ll1$, there is a sufficiently small $\varepsilon^\ast>0$ such that, for $0<\varepsilon<\varepsilon^\ast$, there eixsts a Cantor set $\Pi_\gamma\subset\Pi_0$ with the measure estimate ${\rm meas}\,\Pi_\gamma={\rm meas}\,\Pi_0-O(\gamma)$, and, for each $\sigma\in\Pi_\gamma$, system \eqref{418} has a quasi-periodic invariant torus consisting of quasi-periodic solutions.
\end{thm}
\begin{proof}
System \eqref{418} is a special form of system (1.1) in \cite{LS19} with $n_{11}=n_{21}=0$ (hence, $q_1,q_2$ and $q_6$ can be taken arbitrary values), $n_{12}=n_{22}=2,q_3=q_5=1,q_4=\frac{1}{2}$ and $q_7=\frac{3}{2}$. Taken $q_1=2$ and $q_2=q_6=\frac{1}{2}$, then the hypothesis (H1) on page 4228 in \cite{LS19} is true. As the vector field in \eqref{418} is $C^\infty$ in all variables, the hypothesis (H3) is also true \big(we will take $\alpha=1$ in (H3)\big). Since the eigenvalues $\lambda_1^0,\lambda_2^0$ of $\Omega^0$ satisfy
\begin{equation*}
\lambda_1^0+\lambda_2^0=-2(\rho_{10}+\rho_{20})<0,\qquad \lambda_1^0\lambda_2^0=4\rho_{10}\rho_{20}(1-\delta_1\delta_2)<0,
\end{equation*}
the hypothesis (H2) and these conditions of Remark 2.2 in \cite{LS19} are true. By Theorems 2.2-2.4 and Remark 2.2 in \cite{LS19}, the proof of Theorem \ref{thm5} only remains proving the following non-degeneracy condition
\begin{equation}\label{420}
\det\frac{\partial\omega^1(\sigma)}{\partial\sigma}\neq0,\qquad \forall\sigma\in\Pi_0.
\end{equation}
\begin{proof}[Proof of \eqref{420}]
Let
\begin{equation*}
C=(c_{ij})_{2\times2},\qquad E=(e_{ij})_{2\times2},
\end{equation*}
where $c_{ij},e_{ij}(1\leqslant i,j\leqslant2)$ are defined in Appendix D. Then we have
\begin{equation*}
\dbinom{\tilde\omega_1(\varepsilon\sigma)}{\tilde\omega_2(\varepsilon\sigma)}=\varepsilon EC^{-1}\sigma+O(\varepsilon^2).
\end{equation*}
From \eqref{411}, \eqref{414+} and \eqref{417+}, one can get
\begin{align*}
&s_{11}\rho_{10}+s_{12}\rho_{20}=-\frac{{\rm Im}(d_1e^{-{\rm i}\tau_0\omega_{10}})}{{\rm Re}(d_1e^{-{\rm i}\tau_0\omega_{10}})}\sigma_1+O(\varepsilon)=:d_{11}\sigma_1+O(\varepsilon),\\
&s_{21}\rho_{10}+s_{22}\rho_{20}=-\frac{{\rm Im}(d_2e^{-{\rm i}\tau_0\omega_{20}})}{{\rm Re}(d_2e^{-{\rm i}\tau_0\omega_{20}})}\sigma_2+O(\varepsilon)=:d_{22}\sigma_2+O(\varepsilon).
\end{align*}
Thus,
\begin{equation*}
\omega^1(\sigma)=EC^{-1}\sigma+D_1\sigma
\end{equation*}
and
\begin{equation*}
\frac{\partial\omega^1(\sigma)}{\partial\sigma}=(E+D_1C)C^{-1},
\end{equation*}
where $D_1={\rm diag}(d_{11},d_{22})$.

By using these equalities and expressions given in Appendixes A, D and E, we obtain
\begin{equation*}
\det(E+D_1C)=\frac{2\tau_0(\omega_{10}+\omega_{20})^3(\omega_{10}-\omega_{20})^3(d_{1r}+d_{2r})}{(3\omega_{10}^2+\omega_{20}^2)(\omega_{10}^2+3\omega_{20}^2)\alpha_{30}^2\beta_1^2\beta_2^2d_{1r}d_{2r}}\Delta_2,
\end{equation*}
where
\begin{align*}
d_{1r}=&{\rm Re}(d_1e^{-{\rm i}\tau_0\omega_{10}})\qquad d_{2r}={\rm Re}(d_2e^{-{\rm i}\tau_0\omega_{20}}),\\
\Delta_2=&\alpha_2\left(1+(1+\tau_0\alpha_2)^2+(1+\tau_0\alpha_2)(2+\alpha_{30}^2+\alpha_2^2)\right)^2+\tau_0^4\alpha_2(\alpha_{30}^2+\alpha_2^2)(\alpha_{30}^2-\alpha_2^2)\\
&\hspace{2cm}+\tau_0^2\alpha_2(\alpha_{30}^2-\alpha_2^2)(4+\alpha_{30}^2-\alpha_2^2)(2\tau_0\alpha_2+2+\alpha_{30}^2-\alpha_2^2).
\end{align*}
Therefore,
\begin{equation*}
\det\frac{\partial\omega^1(\sigma)}{\partial\sigma}=\frac{2(d_{1r}+d_{2r})}{\tau_0\alpha_{30}\beta_1\beta_2d_{1r}d_{2r}}\frac{\Delta_2}{\Delta_1}>0.
\end{equation*}
Here, we used \eqref{46*} and facts $\Delta_1,\Delta_2>0$ and $d_{1r},d_{2r}>0$ (see Appendix D), $\Delta_1$ and $\beta_j(j=1,2)$ are defined in \eqref{44}.
\end{proof}
Hence, the proof of Theorem \ref{thm5} is complete.
\end{proof}
By Theorem \ref{thm5}, we obtain the existence of quasi-periodic solutions of system \eqref{22}.
\begin{coro}
For fixed $\alpha_2$, $(\tau_0,\alpha_{30})$ satisfying the conditions in {\bf Case \Rmnum{2}} and \eqref{48}, system \eqref{22} possesses a quasi-periodic solution for most parameter values $\sigma$ (in the sense of Lebesgue measure) of the region \Rmnum{4} depicted in Figure 4.
\end{coro}

\section*{Appendix}
\subsection*{Appendix A: Equalities of $\omega_{10}$ and $\omega_{20}$}
\begin{align*}
&\omega_{10}^2=\frac{1}{2}\left(2+\alpha_{30}^2-\alpha_2^2+\sqrt{(\alpha_{30}^2-\alpha_2^2)(4+\alpha_{30}^2-\alpha_2^2)}\right),\\
&\omega_{20}^2=\frac{1}{2}\left(2+\alpha_{30}^2-\alpha_2^2-\sqrt{(\alpha_{30}^2-\alpha_2^2)(4+\alpha_{30}^2-\alpha_2^2)}\right),\\
&\omega_{10}\omega_{20}=1,\qquad \omega_{10}+\omega_{20}=\sqrt{4+\alpha_{30}^2-\alpha_2^2},\qquad \omega_{10}-\omega_{20}=\sqrt{\alpha_{30}^2-\alpha_2^2},\\
&\omega_{10}^2+\omega_{20}^2=2+\alpha_{30}^2-\alpha_2^2,\qquad \omega_{10}^2-\omega_{20}^2=\sqrt{(\alpha_{30}^2-\alpha_2^2)(4+\alpha_{30}^2-\alpha_2^2)},\\
&(\omega_{10}^2-1)(1+\omega_{20}^2+\tau_0\alpha_2\omega_{20}^2)=\omega_{10}^2-\omega_{20}^2+\tau_0\alpha_2(1-\omega_{20}^2),\\
&(\omega_{20}^2-1)(1+\omega_{10}^2+\tau_0\alpha_2\omega_{10}^2)=\omega_{20}^2-\omega_{10}^2+\tau_0\alpha_2(1-\omega_{10}^2),\\
&(1+\omega_{10}^2+\tau_0\alpha_2\omega_{10}^2)(1+\omega_{20}^2+\tau_0\alpha_2\omega_{20}^2)=1+(1+\tau_0\alpha_2)^2+(1+\tau_0\alpha_2)(\omega_{10}^2+\omega_{20}^2).
\end{align*}

\subsection*{Appendix B: The coefficients of $F(\Phi z,0)$ in \eqref{43}}
\begin{align*}
&F_{3000}=\overline{F_{0300}}=a_1\dbinom{0}{-e^{-3{\rm i}\tau_0\omega_{10}}},\qquad F_{2100}=\overline{F_{1200}}=a_1\dbinom{0}{3e^{-{\rm i}\tau_0\omega_{10}}},\\
&F_{2010}=\overline{F_{0201}}=a_2\dbinom{0}{-e^{-{\rm i}\tau_0(2\omega_{10}+\omega_{20})}},\qquad F_{2001}=\overline{F_{0210}}=a_2\dbinom{0}{e^{-{\rm i}\tau_0(2\omega_{10}-\omega_{20})}},\\
&F_{1110}=\overline{F_{1101}}=a_2\dbinom{0}{2e^{-{\rm i}\tau_0\omega_{20}}},\qquad F_{1020}=\overline{F_{0102}}=a_3\dbinom{0}{-e^{-{\rm i}\tau_0(\omega_{10}+2\omega_{20})}},\\
&F_{1011}=\overline{F_{0111}}=a_3\dbinom{0}{2e^{-{\rm i}\tau_0\omega_{10}}},\qquad F_{1002}=\overline{F_{0120}}=a_3\dbinom{0}{-e^{-{\rm i}\tau_0(\omega_{10}-2\omega_{20})}},\\
&F_{0030}=\overline{F_{0003}}=a_4\dbinom{0}{-e^{-3{\rm i}\tau_0\omega_{20}}},\qquad F_{0021}=\overline{F_{0012}}=a_4\dbinom{0}{3e^{-{\rm i}\tau_0\omega_{20}}},
\end{align*}
where
\begin{align*}
&a_1={\rm i}\tau_0\alpha_4\omega_{10}^3,\quad a_4={\rm i}\tau_0\alpha_4\omega_{20}^3,\\
&a_2={\rm i}3\tau_0\alpha_4\omega_{10}^2\omega_{20}={\rm i}3\tau_0\alpha_4\omega_{10},\\
&a_3={\rm i}3\tau_0\alpha_4\omega_{10}\omega_{20}^2={\rm i}3\tau_0\alpha_4\omega_{20}.
\end{align*}

\subsection*{Appendix C: The expression of $J_1(\mu)$ in \eqref{46}}
\begin{equation*}
J_1(\mu)=(J_{ij})_{4\times4}+O(|\mu|^2)
\end{equation*}
with
\begin{align*}
&J_{11}=\overline{J_{22}}={\rm i}d_1\omega_{10}(1+\omega_{10}^2)\mu_1+d_1\tau_0\omega_{10}^2\mu_2e^{-{\rm i}\tau_0\omega_{10}},\\
&J_{12}=\overline{J_{21}}=-{\rm i}d_1\omega_{10}(1-\omega_{10}^2)\mu_1-d_1\tau_0\omega_{10}^2\mu_2e^{{\rm i}\tau_0\omega_{10}},\\
&J_{13}=\overline{J_{24}}={\rm i}2d_1\omega_{20}\mu_1+d_1\tau_0\mu_2e^{-{\rm i}\tau_0\omega_{20}},\qquad J_{14}=\overline{J_{23}}=-d_1\tau_0\mu_2e^{{\rm i}\tau_0\omega_{20}}\\
&J_{31}=\overline{J_{42}}={\rm i}2d_2\omega_{10}\mu_1+d_2\tau_0\mu_2e^{-{\rm i}\tau_0\omega_{10}},\qquad J_{32}=\overline{J_{41}}=-d_2\tau_0\mu_2e^{{\rm i}\tau_0\omega_{10}}\\
&J_{33}=\overline{J_{44}}={\rm i}d_2\omega_{20}(1+\omega_{20}^2)\mu_1+d_2\tau_0\omega_{20}^2\mu_2e^{-{\rm i}\tau_0\omega_{20}},\\
&J_{34}=\overline{J_{43}}=-{\rm i}d_2\omega_{20}(1-\omega_{20}^2)\mu_1-d_2\tau_0\omega_{20}^2\mu_2e^{{\rm i}\tau_0\omega_{20}}.
\end{align*}
Here, we used the fact that $\omega_{10}\omega_{20}=1$.

\subsection*{Appendix D: The expressions of $\sigma_j(\mu)$ and $\omega_j(\mu)$, $j=1,2$}
Let $\lambda_j=\sigma_j(\mu)+{\rm i}\omega_j(\mu),j=1,2$ be the eigenvalues of the matrix $\tilde J(\mu)$ in \eqref{46}. Then
\begin{equation*}
\sigma_1(\mu)=c_{11}\mu_1+c_{12}\mu_2+O(|\mu|^2),\qquad \sigma_2(\mu)=c_{21}\mu_1+c_{22}\mu_2+O(|\mu|^2)
\end{equation*}
with
\begin{align*}
&c_{11}=-\frac{1}{3\omega_{10}^2+\omega_{20}^2}\big(\omega_{10}(1+\omega_{10}^2)(\omega_{10}^2+\omega_{20}^2){\rm Im}d_1+2(\omega_{10}+\omega_{20}){\rm Im}d_2\big),\\
&c_{12}=\frac{\tau_0}{3\omega_{10}^2+\omega_{20}^2}\big(2{\rm Re}(d_2e^{-{\rm i}\tau_0\omega_{20}})+(1+\omega_{10}^4){\rm Re}(d_1e^{-{\rm i}\tau_0\omega_{10}})\big),\\
&c_{21}=-\frac{1}{\omega_{10}^2+3\omega_{20}^2}\big(\omega_{20}(1+\omega_{20}^2)(\omega_{10}^2+\omega_{20}^2){\rm Im}d_2+2(\omega_{10}+\omega_{20}){\rm Im}d_1\big),\\
&c_{22}=\frac{\tau_0}{\omega_{10}^2+3\omega_{20}^2}\big(2{\rm Re}(d_1e^{-{\rm i}\tau_0\omega_{10}})+(1+\omega_{20}^4){\rm Re}(d_2e^{-{\rm i}\tau_0\omega_{20}})\big),
\end{align*}
and
\begin{equation*}
\omega_1(\mu)=e_{11}\mu_1+e_{12}\mu_2+O(|\mu|^2),\qquad \omega_2(\mu)=e_{21}\mu_1+e_{22}\mu_2+O(|\mu|^2)
\end{equation*}
with
\begin{align*}
&e_{11}=\frac{1}{3\omega_{10}^2+\omega_{20}^2}\omega_{10}(1+\omega_{10}^2)(\omega_{10}^2-\omega_{20}^2){\rm Re}d_1,\\
&e_{12}=\frac{\tau_0}{3\omega_{10}^2+\omega_{20}^2}\omega_{10}^2(\omega_{10}^2-\omega_{20}^2){\rm Im}(d_1e^{-{\rm i}\tau_0\omega_{10}}),\\
&e_{21}=\frac{1}{\omega_{10}^2+3\omega_{20}^2}\omega_{20}(1+\omega_{20}^2)(\omega_{20}^2-\omega_{10}^2){\rm Re}d_2,\\
&e_{22}=\frac{\tau_0}{\omega_{10}^2+3\omega_{20}^2}\omega_{20}^2(\omega_{20}^2-\omega_{10}^2){\rm Im}(d_2e^{-{\rm i}\tau_0\omega_{20}}),
\end{align*}
where, for $j=1,2$,
\begin{align*}
d_j&=\frac{1}{\beta_j}(1+\omega_{j0}^2+\tau_0\alpha_2\omega_{j0}^2-{\rm i}\tau_0\omega_{j0}(\omega_{j0}^2-1),\\
d_je^{-{\rm i}\tau_0\omega_{j0}}&=\frac{1}{\alpha_{30}\omega_{j0}\beta_j}\big(\alpha_2\omega_{j0}(1+\omega_{j0}^2+\tau_0\alpha_2\omega_{j0}^2)+\tau_0\omega_{j0}(\omega_{j0}^2-1)^2\\
&\hspace{5cm}+{\rm i}(\omega_{j0}^2+1)(\omega_{j0}^2-1)\big),
\end{align*}
in which $\beta_j$ are defined in \eqref{46+}. Here, we used some equalities in Appendix A. Hence, we also have
\begin{align*}
&{\rm Im}d_2{\rm Re}(e^{-{\rm i}\tau_0\omega_{10}})-{\rm Im}d_1{\rm Re}(e^{-{\rm i}\tau_0\omega_{20}})\\
&\hspace{1cm}=\frac{\tau_0(\omega_{10}-\omega_{20})}{\alpha_{30}\beta_1\beta_2}\big(\alpha_2(4+2\tau_0\alpha_2+\alpha_{30}^2-\alpha_2^2)-2\tau_0(\omega_{10}^2-1)(\omega_{20}^2-1)\big).
\end{align*}

\subsection*{Appendix E: The expressions of $a_{jl}(0)$ and $q_{jl}(0)$ in \eqref{410}}
\begin{align*}
&a_{11}(0)=\frac{3\tau_0\alpha_4}{\alpha_{30}\beta_1}\omega_{10}^4\big(\alpha_2(1+\omega_{10}^2+\tau_0\omega_{10}^2\alpha_2)+\tau_0(\omega_{10}^2-1)^2\big),\\
&a_{12}(0)=\frac{6\tau_0\alpha_4}{\alpha_{30}\beta_1}\omega_{10}^2\omega_{20}^2\big(\alpha_2(1+\omega_{10}^2+\tau_0\omega_{10}^2\alpha_2)+\tau_0(\omega_{10}^2-1)^2\big),\\
&a_{21}(0)=\frac{6\tau_0\alpha_4}{\alpha_{30}\beta_2}\omega_{10}^2\omega_{20}^2\big(\alpha_2(1+\omega_{20}^2+\tau_0\omega_{20}^2\alpha_2)+\tau_0(\omega_{20}^2-1)^2\big),\\
&a_{22}(0)=\frac{3\tau_0\alpha_4}{\alpha_{30}\beta_2}\omega_{20}^4\big(\alpha_2(1+\omega_{20}^2+\tau_0\omega_{20}^2\alpha_2)+\tau_0(\omega_{20}^2-1)^2\big),\\
&q_{11}(0)=\frac{3\tau_0\alpha_4}{\alpha_{30}\beta_1}\omega_{10}^3(\omega_{10}^2-1)(\omega_{10}^2+1),\quad q_{12}(0)=\frac{6\tau_0\alpha_4}{\alpha_{30}\beta_1}\omega_{10}\omega_{20}^2(\omega_{10}^2-1)(\omega_{10}^2+1),\\
&q_{21}(0)=\frac{6\tau_0\alpha_4}{\alpha_{30}\beta_2}\omega_{10}^2\omega_{20}(\omega_{20}^2-1)(\omega_{20}^2+1),\quad q_{22}(0)=\frac{3\tau_0\alpha_4}{\alpha_{30}\beta_2}\omega_{20}^3(\omega_{20}^2-1)(\omega_{20}^2+1),
\end{align*}
where $\beta_j,j=1,2$ are defined in \eqref{46+}.


\begin{thebibliography}{9}

\markboth{Authors' Names}{Paper Title}


\bibitem[Andriacchi {\it et al.} (1977)]{AOG77}Andriacchi, T. P., Ogle, J. A., Galante, J. O. [1977] ``Walking speed as a basis for normal and abnormal gait measurements,'' {\it J. Biomech.} {\bf 10}, 261-268.
\bibitem[Belykh {\it et al.} (2022)]{BDB22}Belykh, I., Daley, K., Belykh, V. [2022] ``Pedestrian-induced bridge instability: the role of frequency ratios,'' {\it Radiophys. Quantum Electron.} {\bf 64}, 700-708.
\bibitem[Belykh {\it et al.} (2017)]{BJB17}Belykh, I., Jeter, R., Belykh, V. [2017] ``Foot force models of crowd dynamics on a wobbly bridge,'' {\it Sci. Adv.} {\bf 3}, e1701512.
\bibitem[Bocian {\it et al.} (2016)]{BBR16}Bocian, M., Brownjohn, J. M. W., Racic, V., Hester, D., Quattrone, A., Monnickendam, R. [2016] ``A framework for experimental determination of localized vertical pedestrian forces on full-scale structures using wireless attitude and heading reference systems,'' {\it J. Sound Vib.} {\bf 376}, 217-243.
\bibitem[Bocian {\it et al.} (2012)]{BMB12}Bocian, M., Macdonald, J. H. G., Burn, J. F. [2012] ``Biomechanically inspired modeling of pedestrian-induced forces on laterally oscillating structures,'' {\it J. Sound. Vib.} {\bf 331}, 3914-3929.
\bibitem[Bocian {\it et al.} (2013)]{BMB13}Bocian, M., Macdonald, J. H. G., Burn, J. F. [2013] ``Biomechanically inspired modeling of pedestrian-induced vertical self-excited forces,'' {\it J. Bridg. Eng.} {\bf 18}, 1336-1346.
\bibitem[Bramburger {\it et al.} (2014)]{BDL14}Bramburger, J., Dionne, B., LeBlanc, V. G. [2014] ``Zero-Hopf bifurcation in the van der Pol oscillator with delayed position and velocity feedback,'' {\it Nonl. Dyn.} {\bf 78}, 2959-2973.
\bibitem[Brownjohn {\it et al.} (2018)]{BCB18}Brownjohn, J. M. W., Chen, J., Bocian, M., Racic, V., Shahabpoor, E. [2018] ``Using inertial measurement units to identify medio-lateral ground reaction forces due to walking and swaying,'' {\it J. Sound Vib.} {\bf 426}, 90-110.
\bibitem[Cao {\it et al.} (2023)]{CLCH23}Cao, L., Li, J., Chen, Y. F., Huang, S. [2023] ``Measurement and application of walking models for evaluating floor vibration,'' {\it Structures} {\bf 50}, 561-575.
\bibitem[Dallard {\it et al.} (2001)]{DFF01}Dallard, P., Fitezpatrick, A. J., Flint, A., Le Bourva, S., Low, A., Ridsdill, S. R. M., Willford, M. [2001] ``The london millennium footbridge,'' {\it Struct. Eng.} {\bf 79}, 17-33.
\bibitem[Ding {\it et al.} (2013)]{DJY13}Ding, Y., Jiang, W., Yu, P. [2013] ``Double hopf bifurcation in delayed Van der Pol-Duffing equation,'' {\it Int. J. Bifurc. Chaos} {\bf 23}, 1350014.
\bibitem[Ebrahimpour {\it et al.} (1996)]{EHSP96}Ebrahimpour, A., Hmam, A., Sack, R. L., Patten, W. N. [1996] ``Measuring and modeling dynamic loads imposed by moving crowds,'' {\it J. Struct. Eng.} {\bf 122}, 1448-1474.
\bibitem[Erlicher {\it et al.} (2013)]{ETA13}Erlicher, S., Trovato, A., Argoul, P. [2013] ``A modified hybrid Van der Pol/Rayleigh model for the lateral pedestrian force on a periodically moving floor,'' {\it Mech. Syst. Signal Process.} {\bf 41}, 485-501.
\bibitem[Fujino \& Dionysius (2016)]{FD16}Fujino, Y. \& Dionysius, M. S. [2016] ``A conceptual review of pedestrian-induced lateral vibration and crowd synchronization problem on footbridges,'' {\it J. Bridg. Eng.} {\bf 21}, C4015001 1-12.
\bibitem[Goldsztein (2015)]{Gol15}Goldsztein, G. H. [2015] ``Lateral oscillations of the center of mass of bipeds as they walk--Inverted pendulum model with two degrees of freedom,'' {\it AIP. Adv.} {\bf 5}, 107208-1-10.
\bibitem[Han {\it et al.} (2021)]{HZJZ21}Han, H., Zhou, D., Ji, T., Zhang, J. [2021] ``Modelling of lateral forces generated by pedestrians walking across footbridges,'' {\it Appl. Math. Model.} {\bf 89}, 1775-1791.
\bibitem[Harper (1962)]{Har62}Harper, F. C. [1962] ``The mechanics of walking,'' {\it Res. Appl. Indust.} {\bf 15}, 23-28.
\bibitem[He {\it et al.} (2012)]{HLS12}He, X., Li, C., Shu, Y. [2012] ``Triple-zero bifurcation in van der Pol's oscillator with delayed feedback,'' {\it Commun. Nonl. Sci. Numer. Simul.} {\bf 17}, 5229-5239.
\bibitem[Ingolfsson {\it et al.} (2011)]{IGRJ11}Ingolfsson, E. T., Georgakis, C. T., Ricciardelli, F., Jonsson, J. [2011] ``Experimental identification of pedestrian-induced lateral forces on footbridges,'' {\it J. Sound Vib.} {\bf 330},1265-1284.
\bibitem[Ingolfsson {\it et al.} (2012)]{IGJ12}Ingolfsson, E. T., Georgakis, C. T., Jonsson, J. [2012] ``Pedestrian-induced lateral vibrations of footbridges: a literature review,'' {\it Eng. Struct.} {\bf 45}, 21-52.
\bibitem[Jiang {\it et al.} (2015)]{JZS15}Jiang, H., Zhang, T., Song, Y. [2015] ``Delay-induced double Hopf bifurcations in a system of two delay coupled van der Pol's-Duffing oscillators,'' {\it Int. J. Bifurc. Chaos} {\bf 25}, 1550058.
\bibitem[Jiang \& Wei (2007)]{JW08}Jiang, W. \& Wei, J. [2008] ``Bifurcation analysis in van der Pol's oscillator with delayed feedback,'' {\it J. Comput. Appl. Math.} {\bf 213}, 604-615.
\bibitem[Jiang \& Yuan (2007)]{JY07}Jiang, W. \& Yuan, Y. [2007] ``Bogdanov-Takens singularity in van der Pol's oscillator with delayed feedback,'' {\it Physica D: Nonl. Phen.} {\bf 227}, 149-161.
\bibitem[Kerr \& Bishop (2001)]{KB01}Kerr, S. C. \& Bishop, N. [2001] ``Human induced loading on flexible staircases,'' {\it Eng. Struct.} {\bf 23}, 37-45.
\bibitem[Kumar {\it et al.} (2017)]{KKE16}Kumar, P., Kumar, A., Erlicher, S. [2017] ``A modified hybrid van der Pol-Duffing-Rayleigh oscillator for modelling the lateral walking force on a rigid floor,'' {\it Physica D: Nonl. Phen.} {\bf 358}, 1-14.
\bibitem[Kumar {\it et al.} (2018)]{KKR18}Kumar, P., Kumar, A., Racic, V. [2018] ``Modeling of longitudinal human walking force using self-sustained oscillator,'' {\it Int. J. Struct. Stab. Dy.} {\bf 18}, 1850080.
\bibitem[Kumar {\it et al.} (2018)]{KKRE18}Kumar, P., Kumar, A., Racic, V., Erlicher, S. [2018] ``Modelling vertical human walking forces using self-sustained oscillator,'' {\it Mech. Syst. Signal Pr.} {\bf 99}, 345-363.
\bibitem[Kuznetsov (1998)]{Kuz98}Kuznetsov, Y. A. [1998] {\it Elements of Applied Bifurcation Theory} (Springer).
\bibitem[Li \& Shang (2019)]{LS19}Li, X. \& Shang, Z. [2019] ``On the existence of invariant tori in non-conservative dynamical systems with degeneracy and finite differentiability,'' {\it Discr. Cont. Dyn. Syst.} {\bf 39}, 4225-4257.
\bibitem[Li \& Yu (Preprint)]{LY24}Li, X. \& Yu, B. ``The existence of quasi-periodic invariant tori and double Hopf bifurcation of van der Pol's oscillator with delayed feedback,'' Preprint.
\bibitem[Ma {\it et al.} (2008)]{MLF08}Ma, S., Lu, Q., Feng, Z. [2008] ``Double Hopf bifurcation for van der Pol-Duffing oscillator with parametric delay feedback control,'' {\it J. Math. Anal. Appl.} {\bf 338}, 993-1007.
\bibitem[Macdonald (2008)]{Ma08}Macdonald, J. H. G. [2008] ``Lateral excitation of bridges by balancing pedestrians,'' {\it Proc. R. Soc. A-Math. Phys. Eng. Sci.} {\bf 465}, 1055-1073.
\bibitem[Mohammed \& Pavic (2021)]{MP21}Mohammed, A. \& Pavic, A. [2021] ``Human-structure dynamic interaction between building floors and walking occupants in vertical direction,'' {\it Mech. Syst. Signal. Pr.} {\bf 147}, 107036.
\bibitem[Nakamura (2004)]{Nak04}Nakamura, S. [2004] ``Model for lateral excitation of footbridges by synchronous walking,'' {\it J. Struct. Eng. ASCE} {\bf 130}, 32-37.
\bibitem[Nakamura \& Kawasaki (2009)]{NK09}Nakamura, S. \& Kawasaki, T. [2009] ``A method for predicting the lateral girder response of footbridges induced by pedestrians,'' {\it J. Constr. Steel Res.} {\bf 65}, 1705-1711.
\bibitem[Racic \& Brownjohn (2011)]{RB11}Racic, V. \& Brownjohn, J. M. W. [2011] ``Stochastic model of near-periodic vertical loads due to humans walking,'' {\it Adv. Eng. Inf.} {\bf 25}, 259-275.
\bibitem[Racic {\it et al.} (2009)]{RPB09}Racic, V., Pavic, A., Brownjohn, J. M. W. [2009] ``Experimental identification and analytical modeling of human walking forces: literature review,'' {\it J. Sound. Vib.} {\bf 326}, 1-49.
\bibitem[Rainer {\it et al.} (1988)]{RPA88}Rainer, J. H., Pernica, G., Allen, D. E. [1988] ``Dynamic loading and response of footbridges,'' {\it Can. J. Civil. Eng.} {\bf 15}, 66-71.
\bibitem[Ricciardelli \& Demartino (2016)]{RD16}Ricciardelli, F. \& Demartino, C. [2016] ``Design of footbridges against pedestrian-induced vibrations,'' {\it ASCE J. Bridge Eng.} {\bf 21}, C4015003.
\bibitem[Ricciardelli \& Pizzimenti (2007)]{RP07}Ricciardelli, F. \& Pizzimenti, A. D. [2007] ``Lateral Walking-Induced Forces on Footbridges,'' {\it J. Bridg. Eng.} {\bf 12}, 677-688.
\bibitem[Shahabpoor \& Pavic (2018)]{SP18}Shahabpoor, E. \& Pavic, A. [2018] ``Estimation of vertical walking ground reaction force in real-life environments using single IMU sensor,'' {\it J. Biomech.} {\bf 79}, 181-190.
\bibitem[Shahabpoor {\it et al.} (2016)]{SPR16}Shahabpoor, E., Pavic, A., Racic, V., [2016] ``Interaction between walking humans and structures in vertical direction : A literature review,'' {\it Shock Vib.}, ID3430285: 1-22.
\bibitem[Wang \& Jiang (2010)]{WJ10}Wang, H. \& Jiang, W. [2010] ``Hopf-pitchfork bifurcation in van der Pol's oscillator with nonlinear delayed feedback,'' {\it J. Math. Anal. Appl.} {\bf 638}, 9-18.
\bibitem[Wei \& Jiang (2005)]{WJ05}Wei, J. \& Jiang, W. [2005] ``Stability and bifurcation analysis in van der Pol's oscillator with delayed feedback,'' {\it J. Sound Vib.} {\bf 283}, 801-819.
\bibitem[Wiggins (2003)]{Wig03}Wiggins, S. [2003] {\it Introduction to Applied Nonlinear Dynamical Systems and Chaos,} (Springer-Verlag).
\bibitem[Wheeler (1982)]{Wh82}Wheeler, J. E. [1982] ``Prediction and control of pedestrian induced vibration in footbridges,'' {\it J. Struct. Div.} {\bf 108}, 2045-2065.
\bibitem[Wu \& Wang (2011)]{WW11}Wu, X. \& Wang, L. [2011] ``Zero-Hopf bifurcation for van der Pol's oscillator with delayed feedback,'' {\it J. Comput. Appl. Math.} {\bf 235}, 2586-2602.
\bibitem[Zhang \& Guo (2013)]{ZG13}Zhang, L. \& Guo, S. [2013] ``Hopf bifurcation in delayed van der Pol's oscillators,'' {\it Nonl. Dyn.} {\bf 71}, 555-568.
\bibitem[Zhen {\it et al.} (2016)]{ZXS16}Zhen, B., Xu, J., Song, Z. [2016] ``Lateral periodic vibrations of footbridges under crowd excitation,'' {\it Nonl. Dyn.} {\bf 86}, 1701-1710.
\bibitem[Zhen {\it et al.} (2013)]{ZXX13}Zhen, B., Xie, W. P., Xu, J. [2013] ``Nonlinear analysis for the lateral vibration of footbridges induced by pedestrians,'' {\it J. Bridg. Eng. ASCE} {\bf 18}, 122-130.
\bibitem[Zivanovi {\it et al.} (2005)]{ZPR05}Zivanovi, C. S., Pavic, A., Reynolds, P. [2005] ``Vibration serviceability of footbridges under human-induced excitation: A literature review,'' {\it J. Sound Vib.} {\bf 279}, 1-74.

\end{thebibliography}
\end{document}